\documentclass[11pt,reqno]{amsart}
\usepackage[T1]{fontenc}
\usepackage[utf8]{inputenc}
\usepackage{lmodern}
\usepackage[english]{babel}
\usepackage[a4paper,margin=26mm]{geometry}
\usepackage{amsmath,amssymb,amsthm,mathtools,mathrsfs}
\usepackage{booktabs,longtable,array,enumitem}
\usepackage{microtype}
\usepackage{needspace}
\usepackage[hidelinks,unicode]{hyperref}
\usepackage{bookmark}
\hypersetup{
 pdftitle={Derivations of Leavitt path algebras: normal forms, multiplier extensions and growth},
 pdfauthor={Viktor Lopatkin},
 pdfsubject={Lie algebras of derivations and first Hochschild cohomology}}
\setlist{topsep=4pt,itemsep=2pt,parsep=0pt,leftmargin=*}
\numberwithin{equation}{section}
\newtheorem{theorem}{Theorem}[section]
\newtheorem{lemma}[theorem]{Lemma}
\newtheorem{proposition}[theorem]{Proposition}
\newtheorem{corollary}[theorem]{Corollary}
\newtheorem{maintheorem}{Theorem}

\theoremstyle{definition}

\newtheorem{example}[theorem]{Example}
\theoremstyle{remark}
\newtheorem{remark}[theorem]{Remark}
\newcommand{\kk}{k}
\newcommand{\HH}{\operatorname{HH}}
\newcommand{\Der}{\operatorname{Der}}
\newcommand{\ad}{\operatorname{ad}}
\newcommand{\Inn}{\operatorname{Inn}}
\newcommand{\im}{\operatorname{im}}
\newcommand{\rk}{\operatorname{rank}}
\newcommand{\tr}{\operatorname{tr}}
\newcommand{\supp}{\operatorname{supp}}
\newcommand{\diag}{\operatorname{diag}}
\newcommand{\gl}{\mathfrak{gl}}
\newcommand{\slie}{\mathfrak{sl}}
\newcommand{\cD}{\mathscr D}
\newcommand{\cB}{\mathcal B}
\newcommand{\cP}{\mathcal P}
\newcommand{\cQ}{\mathcal Q}
\newcommand{\cM}{\mathcal M}
\newcommand{\cO}{\mathcal O}

\newcommand{\ZZ}{\mathbb Z}
\newcommand{\NN}{\mathbb N}
\newcommand{\reg}{\mathrm{reg}}
\newcommand{\fin}{\mathrm{fin}}
\newcommand{\rcf}{\mathrm{rcf}}

\title[Derivations of Leavitt path algebras]
{Derivations of Leavitt path algebras:\\
normal forms, multiplier extensions and growth}
\author{Viktor Lopatkin}
\date{November 27, 2026}
\begin{document}
\begin{abstract}
We give explicit coordinates for the derivation Lie algebra of a
row-finite Leavitt path algebra over a commutative coefficient ring.
Every derivation is uniquely the sum of a commutator with a
zero-diagonal, row-and-column-finite corner matrix and a derivation
parametrized by the regular vertex corners. These coordinates identify
ordinary inner derivations, compute the Lie bracket, and yield a short
exact sequence separating multiplier-implemented outer classes from
the remaining path parameters. For finite graphs, we construct an
integral normal form for the first Hochschild cohomology: mixed
relations admit an explicit contraction, and the residual relations
are incidence matrices reduced by spanning forests. The reduction
records inner implementers and commutes with extension of scalars.
Over any field, the degree-zero component is finite-dimensional if
and only if the Leavitt path algebra has Gelfand--Kirillov dimension
at most two. In this case all homogeneous dimensions are uniformly
bounded and eventually periodic. Infinite examples include finitary
matrix algebras and embeddings of products of matrix Lie algebras
modulo their direct sums.
We also give a precise comparison with the author's earlier derivation
formulas and correct several conclusions concerning innerness,
infinite graphs, and extension to graph \(C^*\)-algebras.
\end{abstract}
\maketitle
\tableofcontents
\clearpage

\section{Introduction}\label{intro:scope}

The explicit description of derivations is a first step towards
understanding the Lie algebra of infinitesimal symmetries of a
Leavitt path algebra. To compute the first Hochschild cohomology as
a Lie algebra, one must also determine all inner relations, choose
independent outer representatives, and evaluate their brackets.
For an infinite graph an additional distinction arises: a pointwise
finite family of inner operators can define a derivation which is
not implemented by any single element of the algebra.

The starting point is the author's paper
\emph{Derivations of Leavitt path algebras},
published in the \emph{Journal of Algebra} \cite{Lopatkin}.
That paper treats row-finite graphs over a commutative coefficient
ring and describes derivations by their coefficients in the normal
path basis; it also introduces the closed-path, ghost-path, and
mixed-path operators used below
\cite[Theorem~2.1 and Corollary~2.1]{Lopatkin}.
We retain these explicit constructions, fix a consistent sign
convention, and organize them by vertex corners. The coordinates
make the support conditions and all ordinary inner relations
explicit. For finite graphs they lead to independent outer
representatives, a normal-form procedure, and graph-theoretic
consequences. Several conclusions of the earlier paper require
correction; these are stated precisely in
Subsection~\ref{intro:earlier-work}.

Let \(A=L_\kk(\Gamma)\), where \(\Gamma=(V,E,s,r)\) is row-finite
and \(\kk\) is a nonzero commutative ring with identity. We use
\(\HH^1_\kk(A)=\Der_\kk(A)/\ad(A)\), including when \(A\) has only
local units. This is the full algebraic cochain convention, with
degree-zero cochains \(A\). For an infinite graph it must be
distinguished from Ext in the category of unitary bimodules; see
Proposition~\ref{inf:unitary-ext}. Put
\[
\cP=\prod_{v\in V_{\reg}}vAv,\qquad
\cQ=\prod_{v\in V}vAv,
\]
and define
\[
(\partial a)_v=a_v-\sum_{s(e)=v}ea_{r(e)}e^*.
\]
Let \(\cM\) be the algebra of row-and-column-finite matrices with
\((u,v)\)-entry in \(uAv\), and let \(\cO\) be its zero-diagonal
submodule. The algebra \(A\) consists of matrices of finite total
support, and \(\cM\) identifies with its multiplier algebra.

\begin{maintheorem}[Coordinates and the multiplier extension]
\label{intro:coordinates}
Every derivation of \(A\) is uniquely of the form
\[
D=\ad_H+D_p,\qquad H\in\cO,\quad p\in\cP,
\]
where
\[
D_p(v)=0,\qquad D_p(e)=p_{s(e)}e,\qquad
D_p(e^*)=-e^*p_{s(e)}.
\]
It is an ordinary inner derivation precisely when \(H\) has finite
total support and \(p=\partial a\) for a finitely supported
diagonal family \(a\).
There is an exact sequence of \(\kk\)-Lie algebras
\[
0\longrightarrow \frac{\cM}{A+Z(\cM)}
\xrightarrow{\ \jmath\ }\HH^1_\kk(A)
\xrightarrow{\ \pi\ }\frac{\cP}{\partial\cQ}
\longrightarrow0,
\]
with
\(\jmath(M+A+Z(\cM))=[\ad_M|_A]\) and
\(\pi([\ad_H+D_p])=p+\partial\cQ\).
\end{maintheorem}

The full bracket is computed in Theorem~\ref{gen:fullbrackettheorem};
in normalized coordinates it is
\[
\{p,q\}_v=D_p(q_v)-D_q(p_v)+q_vp_v-p_vq_v.
\]
The exact-sequence mechanism comes from the nested Lie ideals
\(\ad(A)\subseteq\ad(\cM)|_A\subseteq\Der_\kk(A)\).
The graph coordinates identify its terms and maps explicitly.
No splitting as Lie algebras is assumed.
If the graph is finite, \(\cM=A\), and the left-hand term vanishes.
For infinite graphs it can account for every outer class, as in
the incoming-star example of Proposition~\ref{inf:star}.

For a finite graph the remaining problem is to compute the cokernel
of \(\partial\) in independent coordinates. We first delete sources.
On the resulting graph, choose one incoming edge at each vertex and
one distinguished outgoing edge at each regular vertex. The normal
path basis of \cite[Theorem~3.7]{LN} makes the mixed inner relations
triangular with respect to path depth. An explicit contraction removes
them. The residual relations are finite incidence matrices in each
degree, so spanning forests select a basis over every coefficient ring.

\begin{maintheorem}[Integral normal form]\label{intro:normal-form}
For a finite graph, \(\HH^1_\kk(A)\) is a free \(\kk\)-module with
explicit homogeneous path representatives.
There is a normal-form map \(\Pi\) on the normalized parameters with
\(\ker\Pi=\im\partial\), together with an explicit linear certificate
recovering an inner implementer whenever \(\Pi(p)=0\).
The construction terminates on every input and commutes with
extension of scalars from \(\ZZ\).
\end{maintheorem}

Theorem~\ref{fin:main} specifies the indexing sets.
The mixed contraction is proved in Lemma~\ref{fin:mixed-kernel};
the remaining basis and innerness certificate are given in
Subsection~\ref{fin:forest-section}.
All eliminations use coefficients \(1\) and \(-1\), without division
by cycle lengths or integers that might vanish in \(\kk\).
This is useful for both positive characteristic and non-field bases.
The resulting brackets, dimensions, and rational path-counting
functions are given in Sections~\ref{fin:lie-section}
and~\ref{fin:ranks-dimensions-section}.

The normal form also detects a precise growth threshold.
Write \(\mathfrak g_d=\HH^1_\kk(A)_d\), where the degree is induced
by \(\deg e=1\), \(\deg e^*=-1\).

\begin{maintheorem}[A cohomological growth criterion]\label{intro:growth}
Let \(\Gamma\) be a nonempty finite graph and let \(\kk\) be any field.
The following conditions are equivalent:
\begin{enumerate}[label=\textup{(\roman*)}]
\item \(\dim_\kk\mathfrak g_0<\infty\);
\item every \(\mathfrak g_d\) is finite-dimensional;
\item the dimensions \(\dim_\kk\mathfrak g_d\), \(d\in\ZZ\), are
uniformly bounded;
\item distinct simple cycles are disjoint and no directed path joins
one such cycle to another;
\item \(\operatorname{GKdim}A\le2\).
\end{enumerate}
When these conditions fail, \(\dim_\kk\mathfrak g_0=\aleph_0\).
When they hold, the sequence
\((\dim_\kk\mathfrak g_d)_{d\ge0}\) is eventually periodic,
with an explicit period and onset bound.
\end{maintheorem}

The growth formula for \(A\) is due to Alahmadi, Alsulami, Jain,
and Zelmanov \cite{AAJZ}. Theorem~\ref{gk:main} connects that
formula with independent outer derivations: two interacting cycles
produce infinitely many independent degree-zero classes, whereas
isolated cyclic components impose a uniform bound on mixed depth.
Here the Gelfand--Kirillov dimension is that of \(A\), not of a
universal enveloping algebra of \(\HH^1_\kk(A)\).

The examples illustrate different aspects of the construction.
An oriented cycle gives the Witt algebra with a rescaled internal
grading; exits contribute further classes and matrix Lie algebras.
For earlier work on derivations of the algebraic Toeplitz algebra,
see Gerritzen \cite{Gerritzen}.
For the infinite incoming star, the description agrees with known
results on derivations of infinite matrix rings \cite{Bez}.
For \(r\ge2\), parallel layers give an injection
\[
\frac{\prod_{n\ge0}\gl_r(\kk)}
     {\bigoplus_{n\ge0}\gl_r(\kk)}
\ \lhook\joinrel\longrightarrow\ \HH^1_\kk(A),
\]
which need not exhaust \(\HH^1_\kk(A)\).
Proposition~\ref{inf:layers-multipliers} further shows that a layer
family is implemented by a multiplier exactly when its matrices
are scalar from some point onwards. Every layer derivation is
locally inner, so a nonscalar tail distinguishes local innerness
from implementation by a multiplier.

Over a field, the bimodule resolution of Chen, Li, and Wang
\cite[Theorem~4.1]{CLWshort} identifies
\(\operatorname{Ext}^1_{A^e,\mathrm{un}}(A,A)\) with
\(\cP/\partial\cQ\). Thus it computes the right-hand term of our
exact sequence; the ordinary cochain quotient in the middle can
be larger. Proposition~\ref{inf:unitary-ext} proves this comparison,
including the products forced by infinitely many local units.
For the related problem of Lie ideals and Lie derivations of the
underlying Lie algebra of the algebraic Toeplitz algebra, see
\cite{HHToep}. Throughout the present paper, derivations satisfy
the associative Leibniz identity.
Appendix~\ref{aut:section}
records first-order and formal automorphisms and polynomial path
shears, with an explicit comparison to the constructions of
Kuroda--Nam \cite{KNauto} and Nam--Srivastava--Vien \cite{NSVauto}.
The existence of the known bouquet shears and the infinite-matrix
derivation theorem is not claimed as new. Appendix~\ref{analytic:boundary}
explains the domain distinction required for any comparison with
graph \(C^*\)-algebras.

\Needspace{9\baselineskip}
\subsection{Relation to the earlier paper and corrections}
\label{intro:earlier-work}

All references to \cite{Lopatkin} in this comparison refer to the
published journal version. The path operators remain useful, but
for an infinite graph the families in
\cite[Corollary~2.1]{Lopatkin} must be interpreted with explicit
pointwise finiteness conditions, rather than as ordinary finite
linear combinations. Corollary~\ref{gen:pathfamily} supplies these
conditions. An admissible pointwise sum of inner operators need not
be an ordinary inner derivation.

The following conclusions of the earlier paper are corrected here.
\begin{enumerate}[label=\textup{(\arabic*)}]
\item
Finiteness of the set of predecessor paths is sufficient, but not
necessary, for innerness. The necessity direction of
\cite[Corollary~3.7, p.~77]{Lopatkin} fails for the Toeplitz graph:
if \(a\) is a loop at \(v\) and \(b:v\to w\) is an exit to a sink,
choose \(a\) as the distinguished edge. Then \(bb^*\) is a normal
mixed monomial and \(\cD_{bb^*}=\ad_v\), although
\(v,a,a^2,\ldots\) are infinitely many paths ending at \(v\).
See Corollary~\ref{gen:finiteancestors} and
Subsection~\ref{app:toeplitz}.
\item
The vanishing assertion in
\cite[Corollary~3.8, p.~79]{Lopatkin} does not hold for arbitrary
infinite row-finite graphs under the ordinary cochain convention.
The incoming star has no path of length greater than one, yet
\(\HH^1_\kk(A)\ne0\); see Proposition~\ref{inf:star}.
The finite acyclic conclusion remains valid.
\item
The proposed free families for the Toeplitz algebra in
\cite[Section~6, p.~82]{Lopatkin} contain dependent outer classes.
In our normalized notation,
\[
\cD_{v,a^n}-\cD_{v,a^{n+1}a^*}=\ad_{a^n}
\qquad(n\ge1).
\]
Subsection~\ref{app:toeplitz} gives the missing inner relations and
the resulting Witt algebra, with consistent bracket signs.
\item
The unrestricted extension assertion in
\cite[Theorem~4.4, p.~80]{Lopatkin} is false for derivations defined
on the entire graph \(C^*\)-algebra. For one loop,
\(t\,d/dt\) on \(\mathbb C[t,t^{-1}]\) has no everywhere-defined
complex-linear derivation extension to \(C(S^1)\).
Appendix~\ref{analytic:boundary} gives an elementary proof and
distinguishes this issue from densely defined derivations.
\end{enumerate}

The structural theorems of the present paper are proved directly
from the defining relations and the normal path basis; they do not
depend on these earlier conclusions. The sign comparison for the
ghost operators is made explicitly in
Subsection~\ref{gen:path-signs}.

The paper is organized as follows. Sections~\ref{gen:conventions}
through~\ref{inf:extension} establish the general coordinates,
brackets, and multiplier extension. Sections~\ref{fin:section}
through~\ref{gk:section} develop the finite normal form and growth
criterion. Sections~\ref{app:examples} through~\ref{inf:layers}
give finite and infinite examples.

\section{Leavitt path algebras and path coordinates}\label{gen:conventions}

Throughout the general constructions, \(\kk\) is a nonzero commutative
ring with identity and \(\Gamma=(V,E,s,r)\) is a row-finite directed graph.
Write \(V_{\reg}=\{v\in V:s^{-1}(v)\ne\varnothing\}\).
The Leavitt path algebra \(A=L_\kk(\Gamma)\) is generated by the vertices,
edges, and ghost edges, subject to
\begin{align}
uv&=\delta_{u,v}v,&s(e)e&=e=er(e),\notag\\
r(e)e^*&=e^*=e^*s(e),&
e^*f&=\delta_{e,f}r(e),\label{gen:ck}\\
v&=\sum_{s(e)=v}ee^*&& (v\in V_{\reg}).\notag
\end{align}
Paths are concatenated from left to right. Vertices are paths of length
zero; for \(w=e_1\cdots e_m\), set \(w^*=e_m^*\cdots e_1^*\).
The symbol \(*\) denotes the \(\kk\)-linear anti-involution fixing the
coefficients. A closed path may revisit its vertices; a simple cycle
does not, apart from its coincident initial and terminal vertices.

Finite sums of vertices are local units. In particular,
\begin{equation}\label{gen:corners}
A=\bigoplus_{u,v\in V}uAv
\end{equation}
is an algebraic direct sum: each \(a\in A\) has only finitely many
nonzero components \(uav\). Orthogonality makes these components unique.
We do not adjoin an identity when \(V\) is infinite.
All derivations are \(\kk\)-linear, and our conventions are
\[
\ad_h(a)=ha-ah,\qquad [D,F]=D\circ F-F\circ D,\qquad
\HH^1_\kk(A)=\Der_\kk(A)/\ad(A).
\]
Thus an inner derivation in this quotient must be implemented by an
element of \(A\) itself. Equivalently, we use the full algebraic
Hochschild cochain complex with \(C^0_\kk(A,A)=A\).
The distinction from cohomology in the category of unitary
bimodules for infinite graphs is explained in
Subsection~\ref{inf:cohomology-conventions}.

Choose a distinguished outgoing edge \(\vartheta(v)\) at each
\(v\in V_{\reg}\). The normal \(\kk\)-basis \(\cB\) consists of vertices,
positive-length paths and their ghosts, and monomials \(wh^*\), where
\(w,h\) have positive length and common terminal vertex, and their last
edges are either different or equal and not distinguished.
This is the basis of \cite[Theorem~3.7]{LN}. The lengths of \(w\) and
\(h\) are independent. Each basis word belongs to exactly one corner
\(uAv\).

To relate our notation to the pathwise coefficients in \cite{Lopatkin},
write, for a generator \(x\),
\begin{equation}\label{gen:coefficients}
\begin{split}
D(x)&=\sum_{v\in V}\alpha_v(x)v
 +\sum_{w\in\mathfrak P}\bigl(\beta_w(x)w+\gamma_w(x)w^*\bigr)
 +\sum_{wh^*\in\mathfrak M}\rho_{wh^*}(x)wh^*\\
&=\sum_{b\in\cB}\rho_b^D(x)b.
\end{split}
\end{equation}
Here \(\mathfrak P\) denotes the positive-length paths and
\(\mathfrak M\) the mixed normal monomials.
For each fixed \(x\), only finitely many coefficients are nonzero.
There is no assumption of a common finite support for all generators.

\begin{remark}[Grading]\label{gen:grading}
The standard grading is given by
\(\deg v=0,\ \deg e=1,\ \deg e^*=-1\).
For homogeneous \(a\in A_n\), define \(D_j(a)=(D(a))_{n+j}\).
Taking homogeneous components in the Leibniz identity shows that each
\(D_j\) is a derivation of degree \(j\), and \(D=\sum_jD_j\) pointwise.
When the graph is finite this sum has finite support as a family of
derivations: the finite set \(V\cup E\cup E^*\) consists of homogeneous
generators, and each of their images has only finitely many homogeneous
components. A homogeneous derivation vanishing on these generators
vanishes on \(A\).
This argument does not require the spaces \(A_n\) to be finite-dimensional.
For an infinite graph only pointwise finiteness is automatic.
\end{remark}

\section{Normalized derivations and path operators}
\label{gen:normalizedsection}

A derivation is \emph{normalized} if it vanishes on every vertex.
We denote their Lie algebra by \(\Der_V(A)\); the subscript refers to
vertices, not to homogeneous degree zero.

\begin{theorem}[Free edge values]\label{gen:free}
For every family \(x_e\in s(e)Ar(e)\), \(e\in E\), there is a unique
normalized derivation \(\delta_x\) satisfying \(\delta_x(e)=x_e\).
Its ghost-edge values are
\begin{equation}\label{gen:ghost}
\delta_x(e^*)=-\sum_{s(f)=s(e)}e^*x_f f^*.
\end{equation}
There is no restriction on the number of nonzero \(x_e\), provided
each \(x_e\) is an element of \(A\).
\end{theorem}

\begin{proof}
Row-finiteness makes the sum in \eqref{gen:ghost} finite.
The proposed values have the required corner supports, so the vertex
and support relations are preserved. If \(s(e)=s(f)\), then
\[
\delta_x(e^*)f+e^*\delta_x(f)
=-\sum_{s(g)=s(e)}e^*x_g g^*f+e^*x_f=0.
\]
If the sources differ, both summands vanish by their corner supports.
For \(v\in V_{\reg}\),
\[
\delta_x\left(\sum_{s(e)=v}ee^*\right)
=\sum_{s(e)=v}x_ee^*
 -\sum_{s(e)=s(f)=v}ee^*x_ff^*
=0.
\]
Consequently the derivation on the free algebra with these generator
values annihilates all defining relations and descends to \(A\).

Conversely, for normalized \(D\), differentiating
\(e=s(e)er(e)\) gives \(D(e)\in s(e)Ar(e)\).
Differentiating \(e^*f=\delta_{e,f}r(e)\), multiplying on the right by
\(f^*\), and summing over the edges from \(s(e)\), gives
\[
D(e^*)=-\sum_{s(f)=s(e)}e^*D(f)f^*,
\]
since \(D(e^*)s(e)=D(e^*)\). The values on generators determine \(D\).
\end{proof}

Put
\[
\cP=\prod_{v\in V_{\reg}}vAv.
\]
Each coordinate is an element of a corner of \(A\); infinitely many
coordinates may be nonzero.

\begin{proposition}[Corner parameters]\label{prop:normalized-param}
The formulas
\begin{equation}\label{gen:Dp}
D_p(v)=0,\qquad D_p(e)=p_{s(e)}e,\qquad
D_p(e^*)=-e^*p_{s(e)}
\end{equation}
give mutually inverse \(\kk\)-linear maps
\begin{equation}\label{gen:parameters}
\cP\xrightarrow{\ \sim\ }\Der_V(A),\quad p\longmapsto D_p,
\qquad
p_v=\sum_{s(e)=v}D_p(e)e^*.
\end{equation}
\end{proposition}

\begin{proof}
Apply Theorem~\ref{gen:free} to \(x_e=p_{s(e)}e\).
Its ghost-edge formula reduces to
\[
-e^*p_{s(e)}\sum_{s(f)=s(e)}ff^*=-e^*p_{s(e)}.
\]
Conversely, for normalized \(D\), the element
\(p_v=\sum_{s(e)=v}D(e)e^*\) lies in \(vAv\) and satisfies
\[
p_{s(e)}e=\sum_{s(f)=s(e)}D(f)f^*e=D(e).
\]
Thus \(D_p=D\) by Theorem~\ref{gen:free}. Finally,
\(\sum_{s(e)=v}p_vee^*=p_v\), proving the other inverse identity.
\end{proof}

\begin{remark}\label{rem:why-regular}
A parameter at a sink would occur in none of the values
\eqref{gen:Dp}. Such parameters would be invisible, not constrained to
vanish. Restricting the product to \(V_{\reg}\) removes this redundancy.
In particular, \(\cP=0\) if the graph has no edges.
\end{remark}

\subsection{The original path operators and their signs}
\label{gen:path-signs}

For \(z\in vAv\), let \(\cD_{v,z}\) be the derivation whose only
possibly nonzero parameter is \(p_v=z\), when \(v\in V_{\reg}\).
Set \(\cD_{v,z}=0\) when \(v\) is a sink.
If \(w=w_0\bar w\) and \(h=h_0\bar h\) are positive-length paths
with common source \(v\) and common range, then
\begin{equation}\label{gen:mixed}
\cD_{wh^*}:=\cD_{v,wh^*},\qquad
\cD_{wh^*}(h_0)=w\bar h^*,\qquad
\cD_{wh^*}(w_0^*)=-\bar w h^*.
\end{equation}
All other generator values are zero. Indeed,
\(h^*f=\delta_{h_0,f}\bar h^*\) and
\(f^*w=\delta_{f,w_0}\bar w\).
This notation also applies before the monomial is reduced to normal form.

For a positive-length closed path \(c=c_0\bar c\) based at \(v\), set
\begin{align}
\cD_c&=\cD_{v,c},&
\cD_c(f)&=cf,&\cD_c(c_0^*)&=-\bar c,\label{gen:cycle}\\
\cD_{c^*}&=-\cD_{v,c^*},&
\cD_{c^*}(c_0)&=-\bar c^*,&
\cD_{c^*}(f^*)&=(cf)^*.\label{gen:ghostcycle}
\end{align}
Here \(f\) ranges over the edges starting at \(v\); unlisted values
are zero. We use the sign in the second line throughout the present
paper. The normalized operator \(\cD_{v,c^*}\) agrees with the
ghost-cycle operator in
\cite[Corollary~2.1(4), p.~64, and Theorem~3.2, p.~78]{Lopatkin},
whereas the list in \cite[Section~6, p.~84]{Lopatkin} uses its
negative, as does our notation \(\cD_{c^*}\).
More precisely,
\begin{equation}\label{gen:sharp}
D^\sharp(a)=D(a^*)^*,\qquad
(\cD_{v,z})^\sharp=-\cD_{v,z^*},\qquad
\cD_{c^*}=(\cD_c)^\sharp.
\end{equation}
The Leibniz identity for \(D^\sharp\) follows from
\[
D^\sharp(ab)=\bigl(D(b^*)a^*+b^*D(a^*)\bigr)^*
=aD^\sharp(b)+D^\sharp(a)b.
\]
The remaining identities follow by evaluation on edges.
Conjugation of linear endomorphisms by \(a\mapsto a^*\) preserves
composition order. Consequently
\begin{equation}\label{gen:sharp-bracket}
[D,F]^\sharp=[D^\sharp,F^\sharp].
\end{equation}
In particular, \(\sharp\) preserves the Lie bracket; the bracket
reversal stated in \cite[Remark~2.2, p.~66]{Lopatkin} does not hold
for this operation.

The vertex parameter satisfies
\begin{equation}\label{gen:vertexparameter}
\cD_{v,v}(f)=f,\quad \cD_{v,v}(f^*)=-f^*
\quad(s(f)=v),\qquad
\cD_{v,v}=\sum_{s(e)=v}\cD_{ee^*}.
\end{equation}
For a distinguished edge \(e\), linearity in the parameter gives
\begin{equation}\label{gen:pathreduction}
\cD_{v,(\alpha e)(\beta e)^*}
=\cD_{v,\alpha\beta^*}
 -\sum_{\substack{s(f)=s(e)\\f\ne e}}
 \cD_{v,(\alpha f)(\beta f)^*},
\end{equation}
where \(r(\alpha)=r(\beta)=s(e)\) and the displayed parameters lie
in \(vAv\). This is just the Cuntz--Krieger relation multiplied by
\(\alpha\) and \(\beta^*\).

Expanding each \(p_v\) in the normal basis gives a pointwise finite
path-operator expansion of \(D_p\). The unreduced path notation has
relations such as \eqref{gen:vertexparameter} and
\eqref{gen:pathreduction}; independent coordinates use the normal
bases of the corners.

\section{Vertex normalization and all derivations}\label{gen:fullsection}

For an arbitrary derivation, the values \(a_v=D(v)\) satisfy
\begin{equation}\label{gen:vertexdata}
a_v=a_vv+va_v,\qquad a_vu+va_u=0\quad(v\ne u).
\end{equation}
These conditions are also sufficient for prescribed vertex values.

\begin{lemma}[Vertex correction]\label{gen:vertexcorrection}
If \(a_v\in A\) satisfy \eqref{gen:vertexdata}, then
\begin{equation}\label{gen:correction}
\delta_a(z)=a_uz+za_w\qquad(z\in uAw)
\end{equation}
defines a derivation with \(\delta_a(v)=a_v\).
\end{lemma}

\begin{proof}
Multiplying \(a_v=a_vv+va_v\) on both sides by \(v\) gives
\(va_vv=2va_vv\), hence \(va_vv=0\), without division by \(2\).
The algebraic corner decomposition makes \eqref{gen:correction}
a well-defined linear map.
For \(z\in uAw\) and \(t\in vAb\), the two middle terms of
\(\delta_a(z)t+z\delta_a(t)\) are \(z(a_w+a_v)t\).
If \(w=v\), these vanish because \(wa_ww=0\); the other terms
equal \(\delta_a(zt)\). If \(w\ne v\), then \(zt=0\), and the middle
terms vanish by \(a_wv+wa_v=0\).
Finally, \(\delta_a(v)=a_vv+va_v=a_v\).
\end{proof}

\begin{theorem}[Generator coordinates]\label{gen:generators}
Choose \(a_v\in A\) satisfying \eqref{gen:vertexdata} and arbitrary
\(x_e\in s(e)Ar(e)\). There is a unique derivation with values
\begin{align}
D(v)&=a_v,\notag\\
D(e)&=a_{s(e)}e+ea_{r(e)}+x_e,\label{gen:generatorvalues}\\
D(e^*)&=a_{r(e)}e^*+e^*a_{s(e)}
 -\sum_{s(f)=s(e)}e^*x_ff^*.\notag
\end{align}
Every derivation arises in this way, with unique data
\(a_v=D(v)\) and \(x_e=s(e)D(e)r(e)\).
\end{theorem}

\begin{proof}
The sum \(\delta_a+\delta_x\) has the required values.
Conversely, differentiating the vertex relations gives
\eqref{gen:vertexdata}, and differentiating \(e=s(e)er(e)\) gives
\[
D(e)=a_{s(e)}e+ea_{r(e)}+s(e)D(e)r(e).
\]
Thus \(D-\delta_a\) is normalized with edge values \(x_e\).
Theorem~\ref{gen:free} gives uniqueness.
\end{proof}

\subsection{Corner matrices}
\label{mat:section}

Let \(\cM\) be the algebra of matrices \(H=(H_{uv})_{u,v\in V}\)
with \(H_{uv}\in uAv\) and finitely many nonzero blocks in each row
and each column. Multiplication is the usual matrix multiplication.
Every entry of a product is a finite sum; the support of any row
or column of a product is contained in a finite union of finite supports.
Matrices with finite total block support identify with \(A\) by
\eqref{gen:corners}. Moreover, \(HA,AH\subseteq A\) for every
\(H\in\cM\), so \(\ad_H|_A\) is a derivation.

Let \(\cO\) be the submodule of zero-diagonal matrices and put
\begin{equation}\label{gen:matrices}
\cO_\fin=\bigoplus_{u\ne v}uAv,\qquad
\cQ=\prod_{v\in V}vAv,\qquad
\cQ_\fin=\bigoplus_{v\in V}vAv.
\end{equation}
We identify \(\cQ\) with the diagonal matrices. Thus
\(\cM=\cO\oplus\cQ\) as a \(\kk\)-module, and
\(A=\cO_\fin\oplus\cQ_\fin\).
The identity of \(\cM\) is \(\diag(v)_{v\in V}\).

\begin{proposition}\label{mat:finite-equals}\label{mat:infinite-strict}
If \(V\) is finite, then \(\cM=A\), \(\cO=\cO_\fin\), and
\(\cQ=\cQ_\fin\). If \(V\) is infinite, then \(\cM\ne A\).
\end{proposition}
\begin{proof}
For finite \(V\), any matrix has finite total block support, and
\((H_{uv})\mapsto\sum_{u,v}H_{uv}\) is an algebra isomorphism.
For infinite \(V\), the identity \(\diag(v)\) of \(\cM\) has
infinite support and hence is not in \(A\).
\end{proof}

Even for a finite graph, \(\cP\) and \(\cQ\) have different index
sets when sinks are present. For example, if \(u\to v\) is a single
edge, then \(\cP=\kk u\) whereas \(\cQ=\kk u\times\kk v\).

\begin{theorem}[Unique full coordinates]\label{gen:all}
Every \(D\in\Der_\kk(A)\) has a unique expression
\begin{equation}\label{gen:full}
D=\ad_H+D_p,\qquad H\in\cO,\quad p\in\cP.
\end{equation}
For \(u\ne v\),
\begin{equation}\label{gen:H}
H_{uv}=uD(v)v=-uD(u)v.
\end{equation}
Conversely, every pair \((H,p)\in\cO\oplus\cP\) defines a derivation.
\end{theorem}

\begin{proof}
Put \(a_v=D(v)\). The first expression in \eqref{gen:H} proves
column-finiteness and the second proves row-finiteness, since each
individual \(a_v\) has finite corner support.
Using \(va_vv=0\), block multiplication gives
\[
Hv=a_vv,\qquad vH=-va_v,\qquad [H,v]=a_v.
\]
For \(z\in uAw\), it also gives
\([H,z]=a_uz+za_w=\delta_a(z)\).
Consequently \(D-\ad_H\) is normalized, and
Proposition~\ref{prop:normalized-param} gives a unique \(p\).
If a zero-diagonal matrix commutes with all vertices, all its blocks
vanish. This proves uniqueness of \(H\). The converse follows from
\eqref{gen:Dp} and \([\cM,A]\subseteq A\).
\end{proof}

These are coordinates of a derivation, not the diagonal and
off-diagonal parts of an inner implementing element.
In terms of \eqref{gen:coefficients}, the coefficient of
\(b\in\cB\cap uAv\), \(u\ne v\), in \(H_{uv}\) is
\(\rho_b^D(v)=-\rho_b^D(u)\). The remaining parameters are
\[
x_e=D(e)-He+eH=s(e)D(e)r(e),\qquad
p_v=\sum_{s(e)=v}x_ee^*.
\]
Conversely, these formulas, the Leibniz identity, and normal reduction
recover all original pathwise coefficients.

\begin{corollary}[All derivations in path notation]\label{gen:pathfamily}
Let \(\cB_{\mathrm{off}}\) consist of the normal basis words in corners
\(uAv\) with \(u\ne v\), and put \(\cB_v=\cB\cap vAv\).
Every derivation has a unique expression
\begin{equation}\label{gen:pathexpansion}
D=\sum_{b\in\cB_{\mathrm{off}}}\mu_b\ad_b
 +\sum_{v\in V_{\reg}}\sum_{b\in\cB_v}\lambda_{v,b}\cD_{v,b},
\end{equation}
subject to the following conditions:
\begin{enumerate}
\item for each vertex \(v\), only finitely many nonzero \(\mu_b\)
have \(vb=b\) or \(bv=b\);
\item for each fixed \(v\), only finitely many \(\lambda_{v,b}\)
are nonzero.
\end{enumerate}
The sums are pointwise finite on \(A\).
For a closed ghost path, our convention \eqref{gen:ghostcycle} gives
\(\cD_{v,c^*}=-\cD_{c^*}\).
\end{corollary}

\begin{proof}
Expand the blocks of \(H\) and each \(p_v\) in the corresponding normal
bases. Condition (1) is precisely row- and column-finiteness of \(H\);
condition (2) says \(p_v\in A\).
Only finitely many summands act on any generator, and hence on any
finite product or finite linear combination of generators.
Uniqueness follows from Theorem~\ref{gen:all} and the normal basis.
\end{proof}

\section{Inner derivations and transfer relations}\label{gen:innersection}

Define the linear map \(\partial:\cQ\to\cP\) by
\begin{equation}\label{gen:partial}
(\partial a)_v=a_v-\sum_{s(e)=v}e a_{r(e)}e^*.
\end{equation}
The outgoing sum is finite. Direct evaluation on edges gives
\begin{equation}\label{gen:diagonalinner}
D_{\partial a}=\ad_a|_A\qquad(a\in\cQ).
\end{equation}
Both sides vanish on vertices and send \(e\) to
\(a_{s(e)}e-ea_{r(e)}\), so Theorem~\ref{gen:free} applies.
If a vertex receives infinitely many edges, \(\partial a\) can have
infinite support even for \(a\in\cQ_\fin\). Thus
\(\partial\cQ_\fin\) is always regarded as a submodule of the full
product \(\cP\).

\begin{theorem}[Ordinary innerness]\label{gen:inner}
The derivation with coordinates \((H,p)\) is inner with an implementer
in \(A\) if and only if
\begin{equation}\label{gen:innercriterion}
H\in\cO_\fin,\qquad p\in\partial\cQ_\fin.
\end{equation}
Consequently there is a \(\kk\)-module isomorphism
\begin{equation}\label{gen:outermodule}
\HH^1_\kk(A)\simeq
\cO/\cO_\fin\ \oplus\ \cP/\partial\cQ_\fin.
\end{equation}
The assertion in \eqref{gen:outermodule} is a module decomposition;
its Lie bracket is described in the next section.
\end{theorem}

\begin{proof}
For \(b\in A\), write \(b=b_{\mathrm{off}}+\sum_v a_v\),
where \(a=(vbv)_v\in\cQ_\fin\).
By \eqref{gen:diagonalinner}, the unique coordinates of \(\ad_b\)
are \((b_{\mathrm{off}},\partial a)\).
Conversely, if \eqref{gen:innercriterion} holds and \(p=\partial a\)
with \(a\in\cQ_\fin\), then \(H+\sum_v a_v\in A\) implements \(D\).
The off-diagonal and diagonal parts can be chosen independently,
which proves the module decomposition.
\end{proof}

\begin{proposition}[Incoming-edge transfer]\label{gen:transferprop}
For \(z\in vAv\),
\begin{equation}\label{gen:transfer}
\cD_{v,z}-\sum_{r(e)=v}\cD_{s(e),eze^*}=\ad_z.
\end{equation}
The incoming sum is pointwise finite, even when infinitely many
edges enter \(v\).
\end{proposition}

\begin{proof}
Apply \eqref{gen:diagonalinner} to the diagonal family supported at
\(v\) with value \(z\). Each summand
\(\cD_{s(e),eze^*}\) can act nontrivially only on the generators
\(e,e^*\), by \(e^*f=f^*e=0\) for \(f\ne e\).
Thus the sum is finite on every generator and then on every element.
For a sink, \(\cD_{v,z}=0\), in accordance with our convention.
\end{proof}

\begin{corollary}[Finitely many predecessor paths]\label{gen:finiteancestors}
Let \(R_v\) be the set of all finite real paths ending at \(v\),
including \(v\). If \(R_v\) is finite, then, for \(z\in vAv\),
\begin{equation}\label{gen:finiteancestorformula}
\cD_{v,z}=\ad_{\sum_{\alpha\in R_v}\alpha z\alpha^*}.
\end{equation}
\end{corollary}
\begin{proof}
Sum \eqref{gen:transfer} over the parameters
\(\alpha z\alpha^*\in s(\alpha)As(\alpha)\).
Every positive-length path in \(R_v\) is uniquely obtained by
prepending its first edge to a shorter path in \(R_v\).
All path operators cancel except \(\cD_{v,z}\).
\end{proof}

\begin{corollary}[Finite graphs]\label{gen:finitecoker}
If \(V\) is finite, then
\begin{equation}\label{gen:finitecokerformula}
\HH^1_\kk(A)\simeq
\frac{\displaystyle\bigoplus_{v\in V_{\reg}}vAv}
{\displaystyle\partial\left(\bigoplus_{v\in V}vAv\right)}.
\end{equation}
\end{corollary}
\begin{proof}
Here \(\cM=A\), \(\cO=\cO_\fin\), and \(\cQ=\cQ_\fin\);
apply Theorem~\ref{gen:inner}.
\end{proof}

\section{Lie brackets in path and matrix coordinates}
\label{gen:liesections}

\subsection{Normalized brackets}

\begin{proposition}\label{gen:normalizedbracketprop}
Under \eqref{gen:parameters}, the Lie bracket is
\begin{equation}\label{gen:normalizedbracket}
[D_p,D_q]=D_{\{p,q\}},\qquad
\{p,q\}_v=D_p(q_v)-D_q(p_v)+q_vp_v-p_vq_v.
\end{equation}
In free edge coordinates it is
\begin{equation}\label{gen:edgebracket}
[\delta_x,\delta_y]=\delta_{[x,y]},\qquad
[x,y]_e=\delta_x(y_e)-\delta_y(x_e).
\end{equation}
\end{proposition}
\begin{proof}
The commutator of normalized derivations is normalized. On an edge
with source \(v\), the Leibniz identity gives
\[
[D_p,D_q](e)
=\bigl(D_p(q_v)-D_q(p_v)+q_vp_v-p_vq_v\bigr)e.
\]
Uniqueness in Theorem~\ref{gen:free} proves the first identity.
The second is the same computation with free edge values.
\end{proof}

For the mixed operator in \eqref{gen:mixed}, substitution along paths
gives
\begin{equation}\label{gen:substitution}
\begin{split}
\cD_{wh^*}(\mu\nu^*)
={}&\sum_{\mu=\mu'h_0\mu''}
 \mu'w\bar h^*\mu''\nu^*\\
&-\sum_{\nu=\nu'w_0\nu''}
 \mu(\nu'')^*\bar w h^*(\nu')^*.
\end{split}
\end{equation}
Every occurrence of the indicated edge is counted separately.
The second line uses
\(\nu^*=(\nu'')^*w_0^*(\nu')^*\);
the order of ghost factors is essential. Reduce the resulting
expressions to \(\cB\). Closed-path operators are evaluated in the
same way from \eqref{gen:cycle}--\eqref{gen:ghostcycle}.

For arbitrary \(D,F\), including nonzero vertex values, the original
coefficients satisfy
\begin{equation}\label{gen:coefficientbracket}
\rho_a^{[D,F]}(x)=
\sum_{b\in\cB}
\bigl(\rho_b^F(x)\rho_a^D(b)-\rho_b^D(x)\rho_a^F(b)\bigr).
\end{equation}
The coefficients on \(b\) refer to its normal expansion under
\(D\) or \(F\). Each sum is finite for fixed \(x\).

\subsection{The full bracket and extension to corner matrices}

A normalized derivation preserves every corner:
\(D_p(uav)=uD_p(a)v\). It therefore extends entrywise to \(\cM\),
with \(D_p(M)_{uv}=D_p(M_{uv})\).
Zero blocks remain zero, so row- and column-finiteness are preserved.
The Leibniz identity follows entrywise from finite sums.
For \(D=\ad_H+D_p\), this gives an extension
\begin{equation}\label{gen:multiplier-extension}
\widetilde D(M)=[H,M]+\bigl(D_p(M_{uv})\bigr)_{u,v}
\qquad(M\in\cM).
\end{equation}
Both \(HM\) and \(MH\) belong to \(\cM\).

\begin{theorem}[Full Lie bracket]\label{gen:fullbrackettheorem}
For \((H,p),(K,q)\in\cO\oplus\cP\), put
\[
Z=[H,K]+D_p(K)-D_q(H),\qquad z=(Z_{vv})_{v\in V}\in\cQ.
\]
Then the unique coordinates of the commutator are
\begin{equation}\label{gen:fullbracket}
[(H,p),(K,q)]=\bigl(Z_{\mathrm{off}},\{p,q\}+\partial z\bigr),
\qquad Z_{\mathrm{off}}=Z-\diag(z).
\end{equation}
\end{theorem}
\begin{proof}
In \(\cM\),
\[
[\ad_H,\ad_K]=\ad_{[H,K]},\qquad
[D_p,\ad_K]=\ad_{D_p(K)}.
\]
Together with \eqref{gen:normalizedbracket}, this yields
\[
[\ad_H+D_p,\ad_K+D_q]=\ad_Z+D_{\{p,q\}}.
\]
Separate \(Z\) into its off-diagonal part and diagonal and apply
\eqref{gen:diagonalinner}. All operations are defined in \(\cM\).
The Jacobi identity is inherited from the derivation Lie algebra
through the bijection \eqref{gen:full}.
\end{proof}

Since \(H,K\) have zero diagonal, \(z_v=([H,K])_{vv}\).
This diagonal correction explains why \(\cO\) need not be a Lie
subalgebra. Quotienting by the pairs in \eqref{gen:innercriterion}
gives the full bracket on \(\HH^1_\kk(A)\).
We shall also use
\begin{equation}\label{gen:partialideal}
\{p,\partial a\}=\partial(D_p(a))\qquad(a\in\cQ).
\end{equation}
Indeed, \([D_p,\ad_a]=\ad_{D_p(a)}\), and normalized parameters are
unique. Hence \(\partial\cQ\) and \(\partial\cQ_\fin\) are Lie ideals
in \((\cP,\{\ ,\ \})\); for the latter, entrywise differentiation
preserves finite support. In particular,
\eqref{gen:finitecokerformula} is a Lie algebra isomorphism with
the quotient bracket \eqref{gen:normalizedbracket}.

Finally, \((H^*)_{uv}=(H_{vu})^*\) extends the anti-involution to
\(\cM\), and
\begin{equation}\label{gen:matrixsharp}
(\ad_H)^\sharp=-\ad_{H^*}.
\end{equation}
This follows from
\((Ha^*-a^*H)^*=aH^*-H^*a\) and agrees with \eqref{gen:sharp}.

\section{Multipliers and ordinary outer derivations}
\label{inf:extension}

Let $\Gamma$ be row-finite, with no restriction on the number of vertices.
Recall that $\cM$ consists of matrices $(M_{uv})_{u,v\in V}$, with
$M_{uv}\in uAv$, having finite support in each row and each column.
Its diagonal subalgebra is $\cQ=\prod_{v\in V}vAv$, and its
zero-diagonal submodule is $\cO$. The ideal $A\subseteq\cM$ consists
of matrices with finite total support.

We first identify $\cM$ with the multiplier algebra $M(A)$. A multiplier
is a pair $(L,R)$ of $\kk$-linear maps $A\to A$ satisfying
\[
L(xy)=L(x)y,\qquad R(xy)=xR(y),\qquad xL(y)=R(x)y.
\]
Its matrix entries are
\[
M_{uv}=uL(v)=R(u)v\in uAv.
\]
The equality follows from compatibility. Since $L(v),R(u)\in A$,
the columns and rows have finite support, respectively. Conversely,
matrix multiplication by any $M\in\cM$ defines such a pair, and
the displayed entries recover it uniquely. Thus $\cM=M(A)$.

\begin{theorem}\label{inf:exact}
There is an exact sequence of $\kk$-Lie algebras
\begin{equation}\label{inf:exact-equation}
0\longrightarrow
\frac{\cM}{A+Z(\cM)}
\xrightarrow{\ \jmath\ }
\HH^1_\kk(A)
\xrightarrow{\ \pi\ }
\frac{\cP}{\partial\cQ}
\longrightarrow0.
\end{equation}
The right-hand term carries the bracket induced by
\eqref{gen:normalizedbracket}. The maps are
\[
\jmath\bigl(M+A+Z(\cM)\bigr)=[\ad_M|_A],
\qquad
\pi\bigl([\ad_H+D_p]\bigr)=p+\partial\cQ.
\]
No Lie algebra splitting is asserted in general.
\end{theorem}
\begin{proof}
By Theorem~\ref{gen:all}, write $D=\ad_H+D_p$.
Formula~\eqref{gen:multiplier-extension} extends $D$ to a derivation
$\widetilde D$ of $\cM$. Indeed, entrywise application of $D_p$
preserves the support of every row and column, and $HM,MH\in\cM$.
For $M\in\cM$ and $x\in A$, the Leibniz rule gives
\[
[D,\ad_M](x)=[\widetilde D(M),x].
\]
Hence the commutators with multipliers form a Lie ideal of
$\Der_\kk(A)$ containing $\ad(A)$.

The identity $[\ad_M,\ad_N]=\ad_{[M,N]}$ shows that $\jmath$
preserves brackets. We compute the kernel of
$M\mapsto[\ad_M|_A]$. If $\ad_M|_A=\ad_b$ for some $b\in A$,
then $C=M-b$ commutes with every element of $A$. For $N\in\cM$
and $x\in A$,
\[
[C,N]x=C(Nx)-N(Cx)=(Nx)C-N(xC)=0.
\]
Taking $x$ to be each vertex shows that every column of $[C,N]$
vanishes. Thus $C\in Z(\cM)$, so the kernel is $A+Z(\cM)$.
The converse inclusion is immediate. This proves that $\jmath$
is well defined and injective.

Adding an ordinary inner derivation $\ad_b$, $b\in A$, to
$\ad_H+D_p$ changes $p$ by $\partial(b_{vv})_{v\in V}$.
Therefore $\pi$ is well defined on $\HH^1_\kk(A)$.
The submodule $\partial\cQ$ is a Lie ideal by
\eqref{gen:partialideal}, and the additional diagonal term in
\eqref{gen:fullbracket} belongs to this ideal. Hence $\pi$ is a
Lie homomorphism. Every $p\in\cP$ defines a derivation $D_p$, with
$\pi([D_p])=p+\partial\cQ$, proving surjectivity.

Finally, if $\pi([D])=0$, write $p=\partial a$ with $a\in\cQ$.
Then \eqref{gen:diagonalinner} gives
\[
D=\ad_H+D_p=\ad_{H+a}|_A,
\]
so $[D]\in\im\jmath$. Conversely, writing $M=H+a\in\cO\oplus\cQ$
gives $\ad_M|_A=\ad_H+D_{\partial a}$ and therefore
$\pi([\ad_M|_A])=0$. This proves exactness at the middle term.
\end{proof}

If $V$ is finite, then $\cM=A$ and $\cQ=\cQ_\fin$, so the left-hand
term vanishes. For an infinite graph the distinction between
$\partial\cQ_\fin$ and $\partial\cQ$ is essential: the former
describes normalized ordinary inner derivations, whereas the latter
describes commutators with arbitrary diagonal multipliers.

A commutator with a multiplier is inner on every prescribed finite
set of elements. Indeed, choose a finite sum of vertices $e$ that is
a two-sided local unit for these elements $x$ and for $Mx,xM$.
Then $eMe\in A$ and $[eMe,x]=[M,x]$ for all the chosen $x$.
This argument does not identify all locally inner derivations with
multiplier commutators; \eqref{inf:exact-equation} uses the latter.

\subsection{Ordinary cochains and unitary bimodule cohomology}
\label{inf:cohomology-conventions}

For an algebra without an identity, the choice of cochain complex
matters already in degree one. Our convention is the full algebraic
Hochschild cochain complex
\[
 C^0_\kk(A,A)=A,\qquad
 C^n_\kk(A,A)=\operatorname{Hom}_\kk(A^{\otimes_\kk n},A)
 \quad(n\ge1),
\]
with its usual differential. Its one-cocycles are the
$\kk$-linear derivations and its one-coboundaries are the ordinary
inner derivations. Thus its first cohomology is precisely the
quotient $\HH^1_\kk(A)=\Der_\kk(A)/\ad(A)$ used here.

There is a different degree-one invariant for algebras with local
units. For the following comparison assume that $\kk$ is a field.
Write $A^e=A\otimes_\kk A^{\mathrm{op}}$, and let
$\operatorname{Ext}^1_{A^e,\mathrm{un}}(A,A)$ be computed in the
category of unitary $A$-bimodules. Here unitary means
$AX=X=XA$; equivalently, each element of $X$ has a two-sided local
unit given by a finite sum of vertices.

\begin{proposition}\label{inf:unitary-ext}
Let $\kk$ be a field and $\Gamma$ a row-finite graph. There are
isomorphisms of $\kk$-vector spaces
\begin{equation}\label{inf:unitary-ext-formula}
 \operatorname{Ext}^1_{A^e,\mathrm{un}}(A,A)
 \simeq \cP/\partial\cQ
 \simeq
 \Der_\kk(A)/\{\ad_M|_A:M\in\cM\}.
\end{equation}
Under the degree-one Lie bracket given by commutators of normalized
derivations, these are also Lie algebra isomorphisms.
In particular, over a field the right-hand term of
\eqref{inf:exact-equation} computes unitary bimodule cohomology,
whereas its middle term retains the ordinary cochain convention.
For a finite graph the two degree-one invariants coincide.
\end{proposition}
\begin{proof}
Set $S=\bigoplus_{v\in V}\kk v$.
In our convention for path multiplication, the projective bimodule
resolution of \cite[Theorem~4.1]{CLWshort} is
\begin{equation}\label{inf:clw-resolution}
 0\longrightarrow
 \bigoplus_{v\in V_{\reg}}Av\otimes_\kk vA
 \xrightarrow{\ d\ }
 \bigoplus_{v\in V}Av\otimes_\kk vA
 \xrightarrow{\ \mu\ }A\longrightarrow0,
\end{equation}
where $\mu$ is multiplication and
\[
 d(v\otimes v)=v\otimes v-\sum_{s(e)=v}e\otimes e^*.
\]
Here $e\otimes e^*$ belongs to the summand indexed by $r(e)$.
The displayed formula is the form of the cited resolution after
reversing its path-multiplication convention.
Each $Av\otimes_\kk vA$ is the principal projective unitary
$A^e$-module $A^e(v\otimes v^{\mathrm{op}})$.

A bimodule map $Av\otimes_\kk vA\to A$ is determined by the
arbitrary value $a_v\in vAv$ of $v\otimes v$:
it sends $xv\otimes vy$ to $xa_vy$.
Since Hom out of a direct sum is a product, applying
$\operatorname{Hom}_{A-A}(-,A)$ to the two projective terms of
\eqref{inf:clw-resolution} gives
\[
 \prod_{v\in V}vAv\xrightarrow{\ d^*\ }
 \prod_{v\in V_{\reg}}vAv.
\]
Evaluation on $d(v\otimes v)$ yields
\[
 (d^*a)_v=a_v-\sum_{s(e)=v}ea_{r(e)}e^*
          =(\partial a)_v.
\]
Its cokernel is therefore the first vector space in
\eqref{inf:unitary-ext-formula}.
The products here are essential when the vertex set is infinite.

For the bracket, use the $S$-relative cochain model associated to the
relative bar resolution, as in \cite[Section~2]{CLWshort}.
Its degree-zero cochains are
$\operatorname{Hom}_{S-S}(S,A)=\prod_vvAv$, and its one-cocycles
are the derivations vanishing on $S$, namely $\Der_V(A)$.
The one-coboundaries are their commutators with diagonal families
in $\cQ$ (up to the harmless sign convention for the cochain
differential). Their normalized parameters are $\partial\cQ$.
The commutator of derivations induces the degree-one bracket,
which under $p\mapsto D_p$ is exactly
\eqref{gen:normalizedbracket}. It descends because
$\{p,\partial a\}=\partial(D_p(a))$.
This specifies the Lie structure used on the Ext group; the
resolution calculation by itself is a vector-space computation.

Finally, Theorem~\ref{gen:all} normalizes an arbitrary derivation
by subtracting a multiplier commutator, and
\eqref{gen:diagonalinner} identifies the remaining normalized
multiplier commutators with $\partial\cQ$.
This proves the second isomorphism and identifies the comparison
map with $\pi$ in \eqref{inf:exact-equation}.
When $V$ is finite, $\cM=A$, so the ordinary and unitary-bimodule
quotients agree.
\end{proof}

\begin{remark}\label{inf:cohomology-star}
The incoming star in Proposition~\ref{inf:star} separates the two
conventions. For $A=M_{\fin}(\NN_0,\kk)$, multiplication gives a
bimodule isomorphism
\[
 AE_{00}\otimes_\kk E_{00}A\xrightarrow{\ \sim\ }A,
 \qquad E_{i0}\otimes E_{0j}\longmapsto E_{ij}.
\]
The left-hand bimodule is projective, so
$\operatorname{Ext}^1_{A^e,\mathrm{un}}(A,A)=0$.
In contrast, Proposition~\ref{inf:star} gives
\[
 \HH^1_\kk(A)\simeq
 \frac{M_{\rcf}(\NN_0,\kk)}
 {M_{\fin}(\NN_0,\kk)+\kk I}\ne0.
\]
For example, a diagonal matrix with entries alternating between
$0$ and $1$ is not a finitary matrix plus a scalar matrix.
Thus a vanishing assertion for unitary bimodule Ext cannot be
transferred to the ordinary quotient used in this paper.
\end{remark}

\section{Finite graphs: independent representatives and normal forms}
\label{fin:section}

In this section $\Gamma=(V,E,s,r)$ is finite and $\kk$ is a nonzero
commutative ring with identity. Put $A=L_\kk(\Gamma)$. We use the
internal grading $\deg v=0$, $\deg e=1$, $\deg e^*=-1$ and the ordinary
Lie bracket $[D,F]=DF-FD$, with no super signs. Since the graph has
finitely many generators and each of their images under a derivation
has finite homogeneous support,
\[
 \Der_\kk(A)=\bigoplus_{d\in\ZZ}\Der_\kk(A)_d,
 \qquad
 \HH^1_\kk(A)=\bigoplus_{d\in\ZZ}\HH^1_\kk(A)_d.
\]
For a set $S$, write $\kk^{(S)}$ for the free module of finite linear
combinations of elements of $S$.

The general parametrization gives
\begin{equation}\label{fin:complex}
 C^0=\bigoplus_{v\in V}vAv,
 \qquad C^1=\bigoplus_{v\in V_{\rm reg}}vAv,
 \qquad
 (\partial h)_v=h_v-\sum_{s(e)=v}eh_{r(e)}e^*.
\end{equation}
Thus a normalized derivation has a parameter $p\in C^1$, and
\begin{equation}\label{fin:parameter}
 D_p(v)=0,\qquad D_p(e)=p_{s(e)}e,\qquad
 D_p(e^*)=-e^*p_{s(e)},\qquad
 \HH^1_\kk(A)\simeq C^1/\im\partial.
\end{equation}
The last assertion is the finite case of Theorem~\ref{gen:inner}.
For $z\in vAv$, let $\varepsilon_v(z)$ and $\iota_v(z)$ denote the
corresponding single-component vectors in $C^1$ and $C^0$, respectively.
We set $\varepsilon_v(z)=0$ when $v$ is a sink and write
$\cD_{v,z}=D_{\varepsilon_v(z)}$. In particular,
\begin{equation}\label{fin:one-inner}
 \cD_{u,z}-\sum_{r(e)=u}\cD_{s(e),eze^*}
 =\ad_z,\qquad z\in uAu.
\end{equation}
At a sink the first term vanishes, whereas the sum over incoming edges
need not vanish. This identity keeps track of both a relation between
outer classes and an implementer for the inner correction.

\subsection{Statement of the normal-form theorem}
\label{fin:main-statement}

Repeatedly delete vertices having no incoming edges, together with all
edges emitted by them, until none remain. Denote the resulting graph
by $\Gamma_\circ$. If nonempty, it has no sources, but it may have
sinks. The following indexing sets are defined on $\Gamma_\circ$.
Choose a special outgoing edge $\vartheta(v)$ at each regular vertex,
as in the normal algebra basis, and an incoming edge $\lambda(v)$ at
each vertex, so that $r(\lambda(v))=v$.

For $n\geq0$, put
\begin{equation}\label{fin:edge-index}
 \mathcal E_n=
 \{(P,e): |P|=n+1,\ s(P)=s(e),\ r(P)=r(e)\}.
\end{equation}
Let $\mathcal H_d$ consist of pairs of paths $(w,h)$ having the same
source and range, with $|w|,|h|\geq2$ and $|w|-|h|=d$, subject to the
following conditions:
\begin{enumerate}
\item their last edges are not both the same special edge
      $\vartheta(v)$;
\item their first edges $f,g$ do not satisfy
      $f=g=\lambda(r(f))$.
\end{enumerate}
The first condition says that $wh^*$ is a normal algebra word. The
second selects the surviving representatives after mixed inner
relations have been eliminated. In
Subsection~\ref{fin:forest-section} we construct a finite auxiliary
graph $\mathscr G_n$ whose edge set is $\mathcal E_n$. Choose a spanning
forest $\mathscr F_n$ in it and put
\begin{equation}\label{fin:cotree-index}
 J_n=\mathcal E_n\setminus E(\mathscr F_n).
\end{equation}

\begin{theorem}\label{fin:main}
For every finite graph $\Gamma$ and every nonzero commutative ring
$\kk$ with identity, the module $\HH^1_\kk(L_\kk(\Gamma))$ is free.
It vanishes if $\Gamma_\circ$ is empty. Otherwise, a basis of its
degree-$d$ component consists of the classes
\begin{align}
 &[\cD_{s(w),wh^*}],
 &&(w,h)\in\mathcal H_d;\label{fin:main-mixed}\\
 &[\cD_{s(e),Pe^*}],
 &&(P,e)\in J_d,\quad d\geq0;\label{fin:main-positive}\\
 &[\cD_{s(e),eP^*}],
 &&(P,e)\in J_{-d},\quad d<0.\label{fin:main-negative}
\end{align}
Parameters on $\Gamma_\circ$ are extended by zero at the deleted
vertices of the original graph.

There is an explicit normal-form map
$\Pi:C^1\to\mathcal V_\Gamma$, where $\mathcal V_\Gamma$ is the free
module on these representatives, such that
\begin{equation}\label{fin:main-kernel}
 \ker\Pi=\im\partial.
\end{equation}
It yields a terminating reduction procedure, an innerness criterion,
and a computation of the Lie bracket in independent coordinates.
All constructions commute with extension of scalars from $\ZZ$ to
$\kk$.
\end{theorem}

We prove the theorem by first deleting sources, then eliminating
mixed relations, and finally reducing the remaining finite relation
matrices by spanning forests. The construction also records an
implementer whenever the normal form vanishes.

\subsection{Deleting sources}
\label{fin:sources-section}

\begin{proposition}\label{fin:sources}
Let $u$ be a vertex of a finite graph $\Gamma$ with no incoming edges,
and let $\Gamma'$ be obtained by deleting $u$ and its outgoing edges.
Then
\[
 \HH^1_\kk(L_\kk(\Gamma))
 \simeq \HH^1_\kk(L_\kk(\Gamma'))
\]
as internally graded Lie algebras. In normalized parameters the
isomorphism deletes the component at $u$; its inverse is induced by
extension by zero.
\end{proposition}
\begin{proof}
A path starting outside $u$ cannot reach $u$. Consequently, for
$v\ne u$ the corner $vL_\kk(\Gamma)v$ is the corresponding corner in
$L_\kk(\Gamma')$. If $u$ is regular, identity
\eqref{fin:one-inner} gives $\cD_{u,z}=\ad_z$ for every $z\in uAu$.
Thus the entire summand at $u$ in $C^1$ lies in $\im\partial$.
Projection onto the other components intertwines $\partial$ with
the smaller graph's boundary operator $\partial'$. Conversely, a
boundary for $\Gamma'$ lifts to a boundary for $\Gamma$, possibly
with an additional component at $u$; that component is already a
boundary. Projection therefore induces an isomorphism of quotients.

If $u$ is also a sink, it is isolated. Its summand $\kk u$ occurs
only in $C^0$ and has zero boundary, so deleting it does not alter
the quotient. All these maps preserve degree. For zero-extended
parameters the bracket has zero component at $u$, and all remaining
components are computed entirely within the surviving corners.
Formula~\eqref{gen:normalizedbracket} therefore shows that the
isomorphism preserves the Lie bracket.
\end{proof}

The deletion process terminates after at most $|V|$ steps. In
particular, a finite acyclic graph is deleted completely, so its
$\HH^1$ vanishes. Sinks reachable from cycles can remain in a
nonempty $\Gamma_\circ$.

\begin{remark}\label{fin:source-convention}
Until Subsection~\ref{fin:ranks-section}, $\Gamma$ denotes the finite
graph without sources obtained by this deletion process. The choices
$\lambda,\vartheta$ and all coordinates below refer to this graph.
\end{remark}

\subsection{Mixed relations and an explicit contraction}
\label{fin:mixed-section}

The depth of a normal mixed word $z=wh^*$ is
$\min(|w|,|h|)$. Vertices, real paths, and ghost paths have depth zero.
For a normal mixed word of positive depth in $uAu$, write
\begin{equation}\label{fin:mixed-row}
 R_z=\partial\iota_u(z)
 =\varepsilon_u(z)-\sum_{r(f)=u}\varepsilon_{s(f)}(fzf^*).
\end{equation}
Each word $fzf^*=(fw)(fh)^*$ is normal, since prefixing both paths
does not change their last edges, and its depth is one greater than
that of $z$. Distinguish the term
\[
 b_z=\varepsilon_{s(\lambda(u))}(\lambda(u)z\lambda(u)^*).
\]
The corresponding reduction is
\begin{equation}\label{fin:mixed-rule}
 b_z\equiv\varepsilon_u(z)-
 \sum_{\substack{r(f)=u\\f\ne\lambda(u)}}
 \varepsilon_{s(f)}(fzf^*)\pmod{\im\partial}.
\end{equation}
Only the first term on the right can again have a distinguished first
pair, and its depth is smaller. The other terms satisfy the second
condition in the definition of $\mathcal H_d$.

Let $T\subset C^0$ be freely spanned by the normal mixed corner words
of positive depth, and let $U\subset C^0$ be spanned by vertices,
positive-length closed real paths, and their ghosts. The normal basis
gives $C^0=T\oplus U$. Let $L\subset C^1$ be spanned by normal corner
words of depth at most one. Let $H_{\rm mix}\subset C^1$ be spanned
by normal words of depth at least two whose first pair is not
distinguished. Its basis is indexed by $\bigsqcup_d\mathcal H_d$.

Every basis vector outside $L\oplus H_{\rm mix}$ is uniquely of the
form $b_z$: delete the common first edge to recover $z$. Define linear
maps
\[
 N:C^1\longrightarrow L\oplus H_{\rm mix},\qquad
 \mathsf H:C^1\longrightarrow T
\]
by setting $N(b)=b$, $\mathsf H(b)=0$ on the basis of
$L\oplus H_{\rm mix}$ and, recursively, setting
\begin{align}
 N(b_z)&=N(\varepsilon_u(z))-
 \sum_{\substack{r(f)=u\\f\ne\lambda(u)}}
 \varepsilon_{s(f)}(fzf^*),\label{fin:N}\\
 \mathsf H(b_z)&=-\iota_u(z)+\mathsf H(\varepsilon_u(z)).
 \label{fin:contraction}
\end{align}
The only recursive call lowers depth, so both maps are well defined
and every evaluation terminates.

\begin{lemma}\label{fin:mixed-kernel}
Regard $N$ as an endomorphism of $C^1$ with image
$L\oplus H_{\rm mix}$. Then
\begin{equation}\label{fin:homotopy}
 N^2=N,\qquad \partial\mathsf H=1-N,\qquad
 \mathsf H\partial|_T=1_T,\qquad N\partial|_T=0.
\end{equation}
In particular, $\partial|_T$ is injective, and
\begin{equation}\label{fin:mixed-splitting}
 \ker N=\partial T,\qquad
 C^1=\partial T\oplus L\oplus H_{\rm mix}.
\end{equation}
\end{lemma}
\begin{proof}
Since $N$ is the identity on its indicated image, $N^2=N$. The
identity $\partial\mathsf H=1-N$ is immediate on
$L\oplus H_{\rm mix}$. On $b_z$, induction on depth gives
\begin{align*}
 \partial\mathsf H(b_z)
 &=-R_z+\partial\mathsf H(\varepsilon_u(z))\\
 &=b_z+\sum_{\substack{r(f)=u\\f\ne\lambda(u)}}
       \varepsilon_{s(f)}(fzf^*)-N(\varepsilon_u(z))
 =b_z-N(b_z).
\end{align*}
The terms indexed by $f\ne\lambda(u)$ are not distinguished, so
\[
 \mathsf H(R_z)
 =\mathsf H(\varepsilon_u(z))-\mathsf H(b_z)=\iota_u(z).
\]
This proves $\mathsf H\partial|_T=1_T$ and injectivity of
$\partial|_T$. Formula~\eqref{fin:N} also gives $N(R_z)=0$, hence
$\partial T\subseteq\ker N$. If $N(p)=0$, then
$p=\partial\mathsf H(p)$, proving the reverse inclusion. The direct
sum follows because $N$ is a projection onto $L\oplus H_{\rm mix}$.
No step uses division in $\kk$.
\end{proof}

\subsection{Low depth and the residual relations}
\label{fin:low-section}

\begin{lemma}\label{fin:low-basis}
For $n\geq0$, the vectors $\varepsilon_{s(e)}(Pe^*)$ indexed by
$(P,e)\in\mathcal E_n$ form a basis of $L_n$. For $n>0$, the vectors
$\varepsilon_{s(e)}(eP^*)$ with the same index set form a basis of
$L_{-n}$.
\end{lemma}
\begin{proof}
All corners in $L\subset C^1$ are based at regular vertices. We first
check membership. If $Pe^*$ is not normal, then $e=\vartheta(v)$,
where $v=s(e)$, and $P=Qe$ for a closed path $Q$ at $v$ of length
$n$; when $n=0$, take $Q=v$. The Cuntz--Krieger relation gives
\begin{equation}\label{fin:low-normalization}
 Pe^*=Qee^*=Q-\sum_{\substack{s(f)=v\\f\ne e}}Qff^*.
\end{equation}
Every mixed term on the right is normal of depth one because $f$ is
not special. Thus normalization takes one step and remains in $L_n$.

Every normal mixed word of depth one and degree $n\geq0$ already
has the form $Pe^*$ in the assertion. A closed real path $c$ of length
$n$ at a regular vertex $v$, including $c=v$ when $n=0$, satisfies
\[
 c=\sum_{s(e)=v}ce\,e^*.
\]
This proves spanning. To prove independence, consider a relation in
one corner and multiply it on the right by an arbitrary edge $f$
emitted by that vertex. Since $e^*f=\delta_{e,f}r(e)$, the relation
becomes
\[
 \sum_{P:(P,f)\in\mathcal E_n}a_{P,f}P=0.
\]
Real paths are part of the normal $\kk$-basis of $A$, so all
$a_{P,f}$ vanish. Doing this for each $f$ and each corner proves
independence over any $\kk$. The anti-involution preserves depth,
reverses degree, and sends $Pe^*$ to $eP^*$; this proves the negative
degree assertion.
\end{proof}

Let $\mathcal C_n$ be the set of closed paths of length $n$ with
specified starting vertex; set $\mathcal C_0=V$. For $c\in\mathcal C_n$
based at $v$, define a row of an integer matrix $M_n$, whose columns
are indexed by $\mathcal E_n$, by
\begin{equation}\label{fin:matrix-row}
 R_c=\sum_{s(f)=v}(cf,f)-\sum_{r(e)=v}(ec,e).
\end{equation}
Repeated columns are combined with their indicated signs. For $c=v$
we use the usual vertex multiplication; the first sum is empty if
$v$ is a sink. Write $\operatorname{Row}_\kk M_n$ for the module
spanned by these rows over $\kk$.

\begin{proposition}\label{fin:split}
For a finite graph without sources and every $d\in\ZZ$, there is an
isomorphism of $\kk$-modules
\begin{equation}\label{fin:module-split}
 \HH^1_\kk(A)_d\simeq
 \kk^{(\mathcal H_d)}\oplus
 \frac{\kk^{\mathcal E_{|d|}}}{\operatorname{Row}_\kk M_{|d|}}.
\end{equation}
There are no further relations among these representatives.
\end{proposition}
\begin{proof}
Since $C^0=T\oplus U$, Lemma~\ref{fin:mixed-kernel} identifies
$C^1/\partial T$ with $L\oplus H_{\rm mix}$. If $c$ is a vertex or
a closed real path, then $\partial\iota_{s(c)}(c)$ already lies in
$L$. Expanding its first term by
$c=\sum_{s(f)=s(c)}cf\,f^*$ gives exactly the row
\eqref{fin:matrix-row} in the basis of Lemma~\ref{fin:low-basis}.
At a sink that first term is absent. Applying $*$ to the boundaries
of positive closed paths gives the same matrices in negative degrees.

Thus $N(\partial U)=\partial U\subset L$, with homogeneous pieces
precisely the row modules of the $M_n$. Quotienting by these introduces
no relation in $H_{\rm mix}$ and accounts for all remaining relations
in $L$.
\end{proof}

The direct sum in~\eqref{fin:module-split} is a module decomposition;
it does not assert that its mixed summand is a Lie ideal.

\subsection{Auxiliary graphs and spanning forests}
\label{fin:forest-section}

Define a finite directed multigraph $\mathscr G_n$ with vertex set
$\mathcal C_n\sqcup\{*\}$ and edge set $\mathcal E_n$. For $(P,e)$,
let its positive endpoint be
\[
 a(P,e)=
 \begin{cases}
 c,&P=ce,\quad c\in\mathcal C_n,\\
 *,&\text{if no such factorization exists},
 \end{cases}
\]
and its negative endpoint be
\[
 b(P,e)=
 \begin{cases}
 c',&P=ec',\quad c'\in\mathcal C_n,\\
 *,&\text{if no such factorization exists}.
 \end{cases}
\]
Each factorization is unique when it exists. Orient the edge from
$a(P,e)$ to $b(P,e)$ and use the incidence convention $+1$ at its
source and $-1$ at its range. Loops are allowed and give zero columns.

\begin{proposition}\label{fin:incidence}
The matrix $M_n$ is the incidence matrix of $\mathscr G_n$ with the
row at $*$ deleted.
\end{proposition}
\begin{proof}
The column $(P,e)$ occurs positively in row $c$ exactly when $P=ce$,
and negatively in row $c'$ exactly when $P=ec'$. There are no other
occurrences in~\eqref{fin:matrix-row}. Equal endpoints cancel, and an
endpoint at $*$ contributes no retained row. These are precisely the
stated incidence conventions.
\end{proof}

\begin{remark}[The degree-zero matrix]\label{fin:M0}
The columns of $M_0$ are indexed by ordered pairs of parallel edges,
and
\[
 (M_0)_{v,(e,e)}=\delta_{v,s(e)}-\delta_{v,r(e)},\qquad
 (M_0)_{v,(f,e)}=0\quad(f\ne e).
\]
Thus $M_0$ consists of the ordinary incidence matrix of $\Gamma$
together with zero columns for distinct parallel edges. The latter
are loops at $*$ in $\mathscr G_0$; loops of $\Gamma$ also give zero
incidence columns. A zero relation column leaves a free coordinate
in the quotient. Over any field, $\rk M_0=|V|-c$, where $c$ is the
number of weak components: a row dependence with coefficients
$(a_v)$ is equivalent to $a_{s(e)}=a_{r(e)}$ for every edge, hence to
constancy on each weak component.
\end{remark}

Choose a spanning tree in each weak component of $\mathscr G_n$,
ignoring loops, and let $\mathscr F_n$ be their union. Root the
component containing $*$ at $*$ and choose any root in each other
component. This completes the definition of $J_n$.

\begin{lemma}\label{fin:forest}
The classes of the column basis vectors indexed by $J_n$ form a basis
of
\begin{equation}\label{fin:forest-quotient}
 \kk^{\mathcal E_n}/\operatorname{Row}_\kk M_n\simeq\kk^{J_n}
\end{equation}
over every commutative ring $\kk$.
\end{lemma}
\begin{proof}
In a component not containing $*$, the sum of all incidence rows is
zero. Deleting the chosen root row therefore leaves the row module
unchanged. In the component containing $*$ the root row is already
absent. The submatrix of the remaining rows and the tree columns is
square with determinant $\pm1$. Indeed, a nontrivial rooted tree has
a leaf other than its root, whose row in this submatrix has exactly
one nonzero entry, equal to $1$ or $-1$. Expanding along that row
reduces to a smaller tree; the empty matrix has determinant $1$.

Invertible integer row operations consequently turn the tree block
into the identity. Ordering tree columns first, the independent
relations have the form $(I\ C)$. Each tree coordinate is uniquely
expressed in terms of the non-tree coordinates. Moreover, a linear
combination of these rows supported only on non-tree columns must
be zero, since its tree coordinates equal its row coefficients.
The non-tree classes therefore both span and are independent. The
integer operations work over every $\kk$.
\end{proof}

\begin{corollary}\label{fin:pi}
There is an explicitly computable linear map
$\Pi:C^1\to\mathcal V_\Gamma$ onto the free module with basis
\eqref{fin:main-mixed}--\eqref{fin:main-negative}, satisfying
$\ker\Pi=\im\partial$. It is the identity on the chosen basis
representatives.
\end{corollary}
\begin{proof}
First delete sources as in Proposition~\ref{fin:sources}. On the
remaining graph normalize algebra words, apply $N$, express $L$ in
the basis of Lemma~\ref{fin:low-basis}, and eliminate tree coordinates
using $(I\ C)$ from Lemma~\ref{fin:forest}. Retain the mixed
coordinates. The resulting composite is $\Pi$.

On a graph without sources, if $\Pi(p)=0$, then
$N(p)\in\partial U$ by Proposition~\ref{fin:split} and
Lemma~\ref{fin:forest}. Choose $u\in U$ with $N(p)=\partial u$.
Then
\[
 p=\partial\mathsf H(p)+N(p)
  =\partial(\mathsf H(p)+u).
\]
Conversely, $N$ kills $\partial T$ and forest reduction kills
$\partial U$, so $\Pi\partial=0$. Source deletion transfers the
kernel assertion back to the original graph. The chosen
representatives are fixed by $N$ and have only non-tree low
coordinates, so $\Pi$ fixes them.

Every input contains finitely many words and degrees. The depth
recursion terminates, and the matrix $M_n$ is finite for every degree
that occurs. Hence this procedure terminates on each input.
\end{proof}

Proposition~\ref{fin:split}, Lemma~\ref{fin:forest}, and
Corollary~\ref{fin:pi} prove the basis assertion of
Theorem~\ref{fin:main}. Recording the elements $\iota_u(z)$ and the
row combinations used in the reduction also produces an inner
implementer for the difference from the normal representative.

\begin{corollary}[A linear certificate for the normal form]
\label{fin:certificate}
Once all choices are fixed, there is a linear map
$\mathcal R:C^1\to C^0$ such that, if
$s:\mathcal V_\Gamma\to C^1$ sends basis indices to their chosen
parameter representatives, then
\begin{equation}\label{fin:certificate-equation}
 \Pi s=1,\qquad 1-s\Pi=\partial\mathcal R.
\end{equation}
Every value of $\mathcal R$ is obtained by finitely many integer
linear operations on the input coefficients. In particular, for
$D=\ad_{h_0}+D_p$ with $D_p$ normalized,
\[
 D=\ad_{\,h_0+\sum_v(\mathcal R(p))_v}+D_{s\Pi(p)}.
\]
If $\Pi(p)=0$, this gives an implementer in $A$.
\end{corollary}
\begin{proof}
For a graph without sources, forest elimination selects independent
rows with an invertible integer tree block. It therefore defines a
linear map $G:L\to U$ such that
\[
 l-s\Pi(l)=\partial G(l).
\]
Concretely, the coefficients of the required combination of vertices
and closed paths are obtained from the tree coordinates of $l$ using
the inverse integer tree matrix. Extend $G$ by zero on $H_{\rm mix}$.
Then
\[
 \mathcal R=\mathsf H+GN
\]
has the required property, because $1-N=\partial\mathsf H$ and
$N-s\Pi=\partial GN$.

Reinsert deleted sources inductively. Suppose $u$ is such a source,
$\pi$ deletes its parameter, and $a$ is the zero extension of
$\mathcal R'(\pi p)$ from the smaller graph. If $u$ is regular, set
\[
 \mathcal R(p)=a+\iota_u\bigl(p_u-(\partial a)_u\bigr).
\]
Since $\partial\iota_u(z)=\varepsilon_u(z)$, this corrects exactly
the missing component. An isolated vertex does not occur in $C^1$
and requires no correction. This proves the assertion by induction;
$\Pi s=1$ follows from the choice of representatives.
\end{proof}

\begin{remark}[Effectivity over the coefficient ring]
\label{fin:effective-base}
Termination means finitely many operations in $\kk$ and integer
linear transformations. An implementation of the zero test also
requires effective representations of elements of $\kk$ and a
decidable equality test in $\kk$. Freeness and
\eqref{fin:certificate-equation} require no such algorithmic
assumption on the ring.
\end{remark}

\begin{corollary}\label{fin:base-change}
For every finite graph $\Gamma$, the natural map
\begin{equation}\label{fin:base-change-map}
 \kk\otimes_{\ZZ}\HH^1_{\ZZ}(L_{\ZZ}(\Gamma))
 \longrightarrow\HH^1_\kk(L_\kk(\Gamma))
\end{equation}
is an isomorphism of graded Lie algebras.
\end{corollary}
\begin{proof}
Make the same choices of $\vartheta$, $\lambda$, and forests over
$\ZZ$ and $\kk$. The normal algebra bases have the same indices.
All mixed and forest reductions use integer coefficients with pivots
$1$ or $-1$. Theorem~\ref{fin:main} therefore gives bases with the
same indices, which the natural scalar-extension map identifies.
The derivation and bracket formulas are also integral, so the map
preserves the Lie bracket. No flatness assumption on the
$\ZZ$-module $\kk$ is needed.
\end{proof}

\section{The Lie bracket and the finite reduction algorithm}
\label{fin:lie-section}

Finite coordinates are computed on $\Gamma_\circ$ and extended by
zero to the original graph when necessary. In formulas involving a
coordinate vector we use its chosen representative in $C^1$.

\begin{proposition}\label{fin:coordinate-bracket}
The bracket in the basis of Theorem~\ref{fin:main} is
\begin{equation}\label{fin:bracket-coordinate}
 [p,q]_{\mathcal V}=\Pi(\{p,q\}),\qquad
 \{p,q\}_v=D_p(q_v)-D_q(p_v)+q_vp_v-p_vq_v.
\end{equation}
\end{proposition}
\begin{proof}
Formula~\eqref{gen:normalizedbracket} gives
$[D_p,D_q]=D_{\{p,q\}}$. Inner derivations form an ideal since
$[D_p,\ad_h]=\ad_{D_p(h)}$. For diagonal $h=(h_v)$, uniqueness of
normalized parameters also gives
\[
 \{p,\partial h\}=\partial\bigl((D_p(h_v))_v\bigr).
\]
Thus the bracket descends to the quotient identified with
$\mathcal V_\Gamma$ by $\Pi$. Alternation and the Jacobi identity
are inherited from the commutator of endomorphisms over any $\kk$.
\end{proof}

For explicit computations, Leibniz's rule gives, for $z\in uAu$,
\begin{equation}\label{fin:path-insertion}
\begin{split}
 \cD_{u,z}(\alpha\beta^*)={}&
 \sum_{\substack{\alpha=\alpha_L\alpha_R\\
        |\alpha_R|>0,\ s(\alpha_R)=u}}
 \alpha_Lz\alpha_R\beta^*\\
 &-\sum_{\substack{\beta=\beta_L\beta_R\\
        |\beta_R|>0,\ s(\beta_R)=u}}
 \alpha\beta_R^*z\beta_L^*.
\end{split}
\end{equation}
Here $\alpha,\beta$ have the same range, and the left factors in the
splittings may have length zero. The first sum selects an edge in
$\alpha$; the second selects an edge in $\beta$ as read in reverse
under $*$. Every Leibniz term occurs once, including repeated
occurrences of the same edge.

Write $\mathscr D_{P,e}=\cD_{s(e),Pe^*}$ when $P$ is a positive-length
path parallel to $e$. If $P,Q$ are parallel to $e,f$, respectively,
then
\begin{equation}\label{fin:substitution}
 [\mathscr D_{P,e},\mathscr D_{Q,f}]
 =\sum_{Q=Q_1eQ_2}\mathscr D_{Q_1PQ_2,f}
  -\sum_{P=P_1fP_2}\mathscr D_{P_1QP_2,e}.
\end{equation}
Indeed, on real edges $\mathscr D_{P,e}$ replaces $e$ by $P$ and
kills the others. The commutator is the indicated difference of
successive substitutions. Equality on real edges implies equality
of normalized derivations by~\eqref{fin:parameter}. Thus
\eqref{fin:substitution} holds before passing to cohomology. The
bracket of two such nonnegative-degree representatives again has
only low representatives, so only forest reduction is required.

The involution on derivations from~\eqref{gen:sharp} is
\begin{equation}\label{fin:star}
 D^\sharp(x)=(D(x^*))^*,\qquad
 (D_p)^\sharp=D_{-p^*},\qquad p^*=(p_v^*)_v.
\end{equation}
It preserves the Lie bracket and reverses degree. Bracket
compatibility follows by conjugating endomorphisms with the linear
involution $*$. The minus sign in the parameter is essential when
transporting bracket formulas to negative degrees.

\subsection{Computing normal forms and inner implementers}
\label{fin:algorithm-section}

Suppose a derivation $D$ is specified by its finite values on graph
generators. Its class can be computed as follows.
\begin{enumerate}
\item Set $h_0=\sum_{v\in V}D(v)v$ and
      $D_0=D-\ad_{h_0}$. The general parametrization gives
      $D_0(v)=0$ and
      $p_v=\sum_{s(e)=v}D_0(e)e^*$.
\item Normalize the algebra words in $p$, delete sources, and apply
      the mixed recursion~\eqref{fin:N}.
\item Express the low part in the basis of
      Lemma~\ref{fin:low-basis} and eliminate the tree columns of
      the required matrices $M_n$. The output is $\Pi(p)$.
\end{enumerate}
By Corollary~\ref{fin:pi}, $D$ is inner exactly when $\Pi(p)=0$.
Keeping the implementers in~\eqref{fin:one-inner} at each reduction
step simultaneously produces an element $h\in A$ with $D=\ad_h$.

The reconstruction across deleted sources can be written explicitly.
Suppose $u_1,\ldots,u_k$ were deleted in that order, and
$(h_v)_{v\in V(\Gamma_\circ)}$ satisfies
$\partial_\circ h=p|_{\Gamma_\circ}$. In reverse order
$i=k,\ldots,1$, set
\begin{equation}\label{fin:source-implementer}
 h_{u_i}=p_{u_i}+\sum_{s(e)=u_i}eh_{r(e)}e^*
\end{equation}
at regular deleted vertices, and set $h_{u_i}=0$ at deleted sinks.
Every edge from $u_i$ ends at a later-deleted vertex or in
$\Gamma_\circ$: otherwise its range would have had an incoming
edge when it was deleted. Thus the right-hand side is already
defined. This enforces $(\partial h)_{u_i}=p_{u_i}$, while edges
from surviving vertices never enter deleted vertices. Consequently
$\partial h=p$ throughout the original graph, and an implementer
for $D$ is $h_0+\sum_vh_v$.

To compute a bracket, take the selected representatives, apply
\eqref{fin:path-insertion} and~\eqref{fin:bracket-coordinate}, and
perform the same reduction. Finite input, decreasing depth, and
finite matrices prove termination. No polynomial complexity bound
is asserted: the number of paths of a prescribed length may grow
exponentially.

\subsection{A matrix Lie subalgebra in degree zero}
\label{fin:matrix-lie-section}

\begin{corollary}\label{fin:matrix-lie}
Let $\Gamma$ be finite without sources, let $b_{uv}$ be the number
of edges $u\to v$, and let $c$ be the number of weak components.
Then $\HH^1_\kk(L_\kk(\Gamma))_0$ contains a Lie subalgebra
\begin{equation}\label{fin:matrix-lie-quotient}
 \mathfrak e_\Gamma\simeq
 \frac{\displaystyle\bigoplus_{b_{uv}>0}
             \mathfrak{gl}_{b_{uv}}(\kk)}
 {\displaystyle\{((t_u-t_v)I_{b_{uv}})_{u,v}:
                   (t_v)\in\kk^V\}}.
\end{equation}
The denominator is central. Over a field,
\begin{equation}\label{fin:matrix-lie-dimension}
 \dim_\kk\mathfrak e_\Gamma
 =\sum_{u,v}b_{uv}^2-|V|+c.
\end{equation}
If the field has characteristic zero and
$q=\#\{(u,v):b_{uv}>0\}$, then
\begin{equation}\label{fin:matrix-lie-split}
 \mathfrak e_\Gamma\simeq
 \bigoplus_{b_{uv}>0}\mathfrak{sl}_{b_{uv}}(\kk)
 \oplus\kk^{q-|V|+c},
\end{equation}
where the last summand is abelian. For an arbitrary finite graph,
apply these formulas to $\Gamma_\circ$ and transport the subalgebra
through the source-deletion isomorphism.
\end{corollary}
\begin{proof}
For parallel edges $e_1,\ldots,e_{b_{uv}}$, the operators
$\cD_{u,e_i e_j^*}$ satisfy the matrix-unit bracket, by
\eqref{fin:substitution}:
\[
 [\cD_{u,e_i e_j^*},\cD_{u,e_k e_l^*}]
 =\delta_{j,k}\cD_{u,e_i e_l^*}
  -\delta_{l,i}\cD_{u,e_k e_j^*}.
\]
Distinct blocks commute, and together they form the whole low part
of degree zero before quotienting. Proposition~\ref{fin:split} says
that the only inner relations in this part come from vertices.
For $h=\sum_vt_vv$,
\[
 (\partial h)_u
 =\sum_{s(e)=u}(t_u-t_{r(e)})ee^*,
\]
which is exactly the central denominator in
\eqref{fin:matrix-lie-quotient}. The absence of further relations
proves that the quotient embeds in $\HH^1_0$.

Over a field, the kernel of
$(t_v)\mapsto(t_u-t_v)_{b_{uv}>0}$ consists of functions constant
on each weak component, so its dimension is $c$. This gives
\eqref{fin:matrix-lie-dimension}. In characteristic zero,
$\mathfrak{gl}_m(\kk)=\mathfrak{sl}_m(\kk)\oplus\kk I_m$.
Factoring the $q$-dimensional central part by the image of dimension
$|V|-c$ proves~\eqref{fin:matrix-lie-split}.
\end{proof}

The quotient~\eqref{fin:matrix-lie-quotient} remains valid in positive
characteristic. The direct-sum description requires extra hypotheses:
if the characteristic divides $m$, then $I_m\in\mathfrak{sl}_m(\kk)$.
Characteristic-independent dimensions do not imply identical Lie
algebra structures in different characteristics.

\section{Ranks and graded dimensions}
\label{fin:ranks-dimensions-section}

\subsection{Ranks of the finite relation matrices}
\label{fin:ranks-section}

The next result applies to any finite graph, including graphs with
sources: $\mathcal C_n$, $\mathcal E_n$, and $M_n$ are defined for all
of them. However, cohomology computations using
\eqref{fin:module-split} require passing first to $\Gamma_\circ$.

Let $B=(b_{uv})$ be the adjacency matrix and $c$ the number of weak
components. For $n>0$, let $\mu_n$ count those weak components that
are exactly simple directed cycles of length dividing $n$, with no
additional edges. An isolated loop counts as a cycle of length one.

\begin{theorem}\label{fin:rank}
For every finite graph, over every field,
\begin{equation}\label{fin:rank-formula}
 \rk M_n=\tr(B^n)-\mu_n\quad(n>0),\qquad
 \rk M_0=|V|-c.
\end{equation}
Every nonzero invariant factor in the Smith normal form of these
integer matrices is $1$.
\end{theorem}
\begin{proof}
For $n=0$, the ordinary incidence columns and the extra zero columns
give rank $|V|-c$, as in Remark~\ref{fin:M0}.

Let $n>0$ and suppose $\sum_{c\in\mathcal C_n}a_cR_c=0$ is a row
dependence over a field. Choose $c=e_1\cdots e_n$ based at $v$ with
$a_c\ne0$. If an edge $f\ne e_1$ leaves $v$, the column $(cf,f)$
has positive contribution $a_c$ and no negative contribution:
$cf$ begins with $e_1$, not $f$. Its positive occurrence is unique,
which is a contradiction. Thus $e_1$ is the only outgoing edge of
$v$. Similarly, if $g\ne e_n$ enters $v$, the column $(gc,g)$ has
only the negative contribution $-a_c$, since $gc$ does not end in
$g$. Therefore $e_n$ is the only incoming edge of $v$.

The column $(ce_1,e_1)$ gives
\[
 a_{e_1e_2\cdots e_n}=a_{e_2\cdots e_ne_1}.
\]
All cyclic rotations of $c$ therefore have the same nonzero
coefficient. Applying the preceding argument at each starting
vertex shows that every visited vertex has exactly one incoming and
one outgoing edge. These vertices form an entire weak component,
which is a simple directed cycle of length $r$ dividing $n$.

Conversely, such a component has exactly $r$ based closed paths of
length $n$. Their rows have one dependence with equal coefficients,
and the adjacent rotation equations show that there are no other
dependences on that component. Distinct components are independent.
Thus the row-dependence space has dimension $\mu_n$. As the number
of rows is $|\mathcal C_n|=\tr(B^n)$, this proves the rank formula
in every characteristic.

The auxiliary-graph construction and its spanning forests do not
require the absence of sources. After dependent root rows are
deleted, the forest minor has full rank and determinant $\pm1$.
The product of the nonzero Smith invariant factors is the positive
greatest common divisor of all full-rank minors, which is therefore
$1$. Each nonzero invariant factor is consequently $1$.
\end{proof}

The number of columns is
\begin{equation}\label{fin:kappa}
 \kappa_n=|\mathcal E_n|
 =\sum_{u,v}b_{uv}(B^{n+1})_{uv}:
\end{equation}
for each edge $u\to v$ one independently chooses a path of length
$n+1$ with the same endpoints. Hence
\begin{equation}\label{fin:low-rank}
 \ell_n=|J_n|=
 \begin{cases}
 \kappa_n-\tr(B^n)+\mu_n,&n>0,\\
 \displaystyle\sum_{u,v}b_{uv}^2-|V|+c,&n=0.
 \end{cases}
\end{equation}
These are the ranks of the low cohomology summands when $B$ is
computed on a graph without sources. For a graph with sources,
first use Proposition~\ref{fin:sources}.

\subsection{Counting mixed representatives}
\label{fin:counts-section}

Assume again that $\Gamma$ is finite without sources; for an original
graph use all the following data on $\Gamma_\circ$. Let $Q$ be the
diagonal matrix with $Q_{vv}=1$ at regular vertices and $Q_{vv}=0$
at sinks. For matrices of the same size write
\[
 \langle X,Y\rangle=\sum_{u,v}X_{uv}Y_{uv}.
\]
All counts below are ordinary integers, regardless of the
characteristic of $\kk$. Let $h_{r,s}$ be the number of mixed basis
pairs with $|w|=r$, $|h|=s$, for $r,s\geq2$.

\begin{theorem}\label{fin:mixed-count}
For a finite graph without sources,
\begin{equation}\label{fin:mixed-count-formula}
\begin{split}
 h_{r,s}={}&\langle B^r,B^s\rangle
 -\langle B^{r-1},B^{s-1}\rangle\\
 &-\langle B^{r-1}Q,B^{s-1}Q\rangle
 +\langle B^{r-2}Q,B^{s-2}Q\rangle.
\end{split}
\end{equation}
In particular, these numbers are independent of $\lambda$ and
$\vartheta$ and satisfy $h_{r,s}=h_{s,r}$.
\end{theorem}
\begin{proof}
There are $\langle B^r,B^s\rangle$ pairs of parallel paths of the
specified lengths. Let $F$ be the subset with distinguished first
pair $f=g=\lambda(r(f))$, and let $Z$ be the subset with forbidden
last pair $a=b=\vartheta(s(a))$.

Deleting the first edge gives a bijection from $F$ to all parallel
path pairs of lengths $r-1,s-1$: a pair starting at $u$ has the unique
allowed prefix $\lambda(u)$. Here the absence of sources is used.
Thus $|F|=\langle B^{r-1},B^{s-1}\rangle$. Deleting the last edge
gives a bijection from $Z$ to parallel path pairs of lengths
$r-1,s-1$ ending at a regular vertex, so
$|Z|=\langle B^{r-1}Q,B^{s-1}Q\rangle$.

Deleting both first and last edges in $F\cap Z$ gives parallel paths
of lengths $r-2,s-2$ ending at a regular vertex. Both deleted edges
are recovered uniquely, including when one or both remaining paths
have length zero. Hence
$|F\cap Z|=\langle B^{r-2}Q,B^{s-2}Q\rangle$.
Inclusion--exclusion proves the formula, which is visibly independent
of the edge choices and symmetric in $r,s$.
\end{proof}

\begin{theorem}\label{fin:dimension}
Let $\kk$ be a field and $\Gamma$ an arbitrary finite graph. If
$\Gamma_\circ$ is empty, all components of $\HH^1_\kk(L_\kk(\Gamma))$
vanish. Otherwise compute $B,Q,\ell_n,h_{r,s}$ on $\Gamma_\circ$.
Then
\begin{equation}\label{fin:dimension-formula}
 \dim_\kk\HH^1_\kk(L_\kk(\Gamma))_d
 =\ell_{|d|}+\sum_{m=2}^{\infty}h_{m+|d|,m}
 \qquad(d\in\ZZ).
\end{equation}
An infinite sum is interpreted as countably infinite dimension.
These dimensions are independent of the field's characteristic and
are equal in degrees $d$ and $-d$.
\end{theorem}
\begin{proof}
Theorem~\ref{fin:main} gives $\ell_{|d|}$ low basis elements. For
$d\geq0$, the mixed basis is the disjoint union of the finite sets
of path pairs of lengths $(m+d,m)$, of cardinalities $h_{m+d,m}$.
For $d<0$, exchange the two paths and use symmetry. A finite graph
has at most countably many finite paths, so an infinite set of
these representatives has cardinality $\aleph_0$. Their indices and
counts are integral and do not depend on the coefficient field.
\end{proof}

If a homogeneous component is infinite dimensional,
\eqref{fin:dimension-formula} is not an ordinary Hilbert series with
finite integer coefficients; it is a statement about basis
cardinalities.

\subsection{A finite criterion using a graph of vertex pairs}
\label{fin:transfer-section}

For a finite graph without sources, form the graph of ordered vertex
pairs with adjacency matrix
\[
 T_{\rm pair}=B\otimes B.
\]
Its vertices are $(x,y)\in V^2$, and an edge
$(x,y)\to(x',y')$ is a pair of edges $x\to x'$, $y\to y'$.

Let $\mathcal A$ consist of edge pairs $(f,g)$ with $s(f)=s(g)$
other than $f=g=\lambda(r(f))$. Let $\mathcal Z$ consist of edge
pairs $(a,b)$ with $r(a)=r(b)$ other than
$a=b=\vartheta(s(a))$. For $n\geq0$, define a nonnegative integer
row vector $\alpha_n$ and column vector $\beta$, indexed by $V^2$,
by
\begin{align}
 \alpha_n(x,y)&=
 \sum_{\substack{(f,g)\in\mathcal A\\r(g)=y}}
 (B^n)_{r(f),x},\label{fin:alpha}\\
 \beta(x,y)&=\#\{(a,b)\in\mathcal Z:s(a)=x,\ s(b)=y\}.
 \label{fin:beta}
\end{align}

\begin{proposition}\label{fin:transfer}
For $m\geq2$ and $n\geq0$,
\begin{equation}\label{fin:transfer-formula}
 h_{m+n,m}=\alpha_n T_{\rm pair}^{m-2}\beta.
\end{equation}
\end{proposition}
\begin{proof}
Split a path $w$ of length $m+n$ into its first edge $f$, the next
$n$ edges, the next $m-2$ edges, and its last edge $a$. Split a path
$h$ of length $m$ into its first edge $g$, the next $m-2$ edges, and
its last edge $b$. The vector $\alpha_n$ counts an admissible first
pair and the extra $n$ steps in $w$. The matrix
$T_{\rm pair}^{m-2}$ counts the simultaneous middle steps, and
$\beta$ counts admissible last pairs with common range. These
successive choices are in bijection with the mixed pairs counted
by $h_{m+n,m}$.
\end{proof}

\begin{corollary}\label{fin:infinite-criterion}
Let $\kk$ be a field, and compute $B,\alpha_n,\beta$ on the nonempty
$\Gamma_\circ$ of a finite graph $\Gamma$. The component
$\HH^1_\kk(L_\kk(\Gamma))_d$ is infinite dimensional if and only if
the pair graph has a directed cycle reachable from
$\supp\alpha_{|d|}$ from which $\supp\beta$ is reachable. If there
is no such cycle, then
\begin{equation}\label{fin:finite-dimension-bound}
 \dim_\kk\HH^1_\kk(L_\kk(\Gamma))_d
 =\ell_{|d|}+
 \sum_{k=0}^{|V(\Gamma_\circ)|^2-1}
 \alpha_{|d|}T_{\rm pair}^{k}\beta.
\end{equation}
\end{corollary}
\begin{proof}
If such a cycle exists, start from $\supp\alpha_{|d|}$, reach the
cycle, traverse it arbitrarily many times, and continue to
$\supp\beta$. This produces paths of unbounded lengths giving
positive integer contributions to~\eqref{fin:transfer-formula},
and therefore infinitely many mixed basis elements in degree $d$.

Conversely, infinite dimension in~\eqref{fin:dimension-formula}
requires nonzero $\alpha_{|d|}T_{\rm pair}^{k}\beta$ for unbounded
$k$, since only finitely many paths occur at each length. A sufficiently
long path in a finite graph repeats a vertex, hence contains a
directed cycle with the required reachability.

If there is no such cycle, the subgraph of vertices lying on paths
from $\supp\alpha_{|d|}$ to $\supp\beta$ is acyclic. Every path in
it has at most $|V(\Gamma_\circ)|^2-1$ edges. All later terms vanish,
which proves the finite sum formula.
\end{proof}

This criterion is decidable by finite reachability and strongly
connected component computations in the pair graph. A one-vertex
strong component contains a cycle only if it has a loop. Cycles
outside paths between the two specified supports have no effect.

\subsection{A generating function for the two path lengths}
\label{fin:generating-section}

\begin{corollary}\label{fin:generating-function}
For a finite graph without sources, put
$R_X=(I-XB)^{-1}$ and $R_Y=(I-YB)^{-1}$ as matrices of formal power
series. Then
\begin{equation}\label{fin:generating-formula}
\begin{split}
 \sum_{r,s\geq2}h_{r,s}X^rY^s
 =X^2Y^2\big(&\langle B^2R_X,B^2R_Y\rangle
             -\langle BR_X,BR_Y\rangle\\
            &-\langle BR_XQ,BR_YQ\rangle
             +\langle R_XQ,R_YQ\rangle\big).
\end{split}
\end{equation}
This is a rational function in $\mathbb Q(X,Y)$ whose expansion at
the origin has nonnegative integer coefficients.
\end{corollary}
\begin{proof}
Use $R_X=\sum_{j\geq0}X^jB^j$ and the analogous expansion for $R_Y$
to sum the four terms of~\eqref{fin:mixed-count-formula} over
$r,s\geq2$. The first term gives
$X^2Y^2\langle B^2R_X,B^2R_Y\rangle$, and the last gives
$X^2Y^2\langle R_XQ,R_YQ\rangle$; the other two are obtained in
the same way. The determinants of $I-XB$ and $I-YB$ have constant
term $1$, so their inverses have rational entries regular at the
origin. Nonnegative integrality follows from the interpretation as
numbers of mixed basis pairs.
\end{proof}

Formula~\eqref{fin:generating-formula} counts the two lengths
separately. Passing to their difference requires diagonal sums as
in~\eqref{fin:dimension-formula}, and those sums may be infinite.

\section{The degree-zero component and growth of the algebra}
\label{gk:section}

Throughout this section, $\kk$ is a field of arbitrary characteristic,
$\Gamma$ is a nonempty finite graph, $A=L_\kk(\Gamma)$, and
$\mathfrak g=\HH^1_\kk(A)=\bigoplus_{d\in\ZZ}\mathfrak g_d$.
A cycle means a simple directed cycle; cyclic changes of its initial
vertex do not produce different cycles.
Finite-dimensionality below refers to the homogeneous components,
not to finite-dimensionality of finitely generated Lie subalgebras.

For a finite-dimensional generating subspace $F\subseteq A$ containing
$1$, the Gelfand--Kirillov dimension is
\[
 \operatorname{GKdim}A
 =\limsup_{n\to\infty}\frac{\log\dim_\kk F^n}{\log n}.
\]
One may take $F$ to be the span of the vertices, edges, and ghost edges.
A nonzero finite-dimensional algebra has GK-dimension zero.
We use the growth formula of Alahmadi, Alsulami, Jain, and Zelmanov
\cite[Theorem~5]{AAJZ}; the connection with $\mathfrak g_0$ will follow
from the normal form established above.

\begin{theorem}[A cohomological criterion for growth of degree at most two]
\label{gk:main}
For a nonempty finite graph $\Gamma$ over any field $\kk$, the following
conditions are equivalent:
\begin{enumerate}[label=\textup{(\roman*)}]
\item $\dim_\kk\mathfrak g_0<\infty$;
\item $\dim_\kk\mathfrak g_d<\infty$ for every $d\in\ZZ$;
\item the numbers $\dim_\kk\mathfrak g_d$, $d\in\ZZ$, have a common
finite bound;
\item distinct simple cycles have no common vertex, and there is no
directed path from one simple cycle to a different one;
\item $\operatorname{GKdim}A\le2$.
\end{enumerate}
If these conditions fail, then
\begin{equation}\label{gk:infinite-zero}
 \dim_\kk\mathfrak g_0=\aleph_0.
\end{equation}
If they hold, the sequence $\dim_\kk\mathfrak g_d$, $d\ge0$, is eventually
periodic. An explicit period and a sufficient starting index are given
in Proposition~\ref{gk:periodicity}.
\end{theorem}

By Proposition~\ref{fin:sources}, deleting sources preserves
$\mathfrak g$ and its grading. It also preserves all cycles and all
paths between cycles: no vertex reachable from a cycle is deleted.
Thus condition~\textup{(iv)} can be checked on either $\Gamma$ or
$\Gamma_\circ$. We do not need invariance of GK-dimension under source
deletion.

\subsection{Independent classes arising from two cycles}
\label{gk:bad-cycles}

\begin{lemma}\label{gk:witnesses}
If condition~\textup{(iv)} of Theorem~\ref{gk:main} fails, then
$\mathcal H_0$ contains infinitely many pairs whose first edges differ
and whose last edges differ. Consequently, the associated classes
$[\cD_{s(w),wh^*}]$ are linearly independent for every choice of
$\lambda,\vartheta$ and over every field.
\end{lemma}
\begin{proof}
All cycles and connecting paths used below survive in $\Gamma_\circ$.
First suppose that two vertex-disjoint cycles $C,D$ are joined by a
directed path. Choose a shortest path $P$ from the vertex set of $C$
to that of $D$, and base the cycles at $s(P),r(P)$.
Write $a=|C|$, $b=|D|$. The first edge of $P$ is not an edge of $C$,
and its last edge is not an edge of $D$, by minimality.
For $k\ge1$, set
\begin{equation}\label{gk:disjoint-witness}
 w_k=C^{bk}P,\qquad h_k=PD^{ak}.
\end{equation}
These paths have the same initial vertex, terminal vertex, and length
$abk+|P|$. Their first edges differ, as do their last edges.
Their lengths are at least two, so $(w_k,h_k)\in\mathcal H_0$.
Different values of $k$ give different lengths.

Next suppose that two distinct simple cycles intersect. Base them at
a common vertex where their first edges differ. Such a vertex exists:
otherwise the successor edge would agree at every common vertex and
the cycles would coincide. Again write $a=|C|$, $b=|D|$, and let $R$
be their longest common suffix:
\[
 C=UR,\qquad D=VR.
\]
The suffix may have length zero. Both $U$ and $V$ are nonempty, since
a simple cycle based at the chosen vertex cannot be a proper suffix
of another such cycle. Their last edges differ by maximality of $R$.
For $k\ge2$, set
\begin{equation}\label{gk:intersecting-witness}
 w_k=C^{bk-1}U,\qquad h_k=D^{ak-1}V.
\end{equation}
These paths have common initial and terminal vertices and equal
length $abk-|R|\ge2$. Each starts with a complete copy of its cycle,
so the first edges differ; the last edges differ by the choice of
$U,V$. They therefore give distinct elements of $\mathcal H_0$.

In both cases, linear independence of the outer classes follows
from Theorem~\ref{fin:main}. Independence of the corresponding words
in $A$ alone would not exclude inner relations.
\end{proof}

\subsection{Bounded mixed depth}
\label{gk:bounded-depth}

\begin{lemma}\label{gk:core-shape}
Suppose $\Gamma_\circ\ne\varnothing$ and condition~\textup{(iv)} holds.
Let $a$ be the number of vertices of $\Gamma_\circ$ lying on no cycle.
Then:
\begin{enumerate}[label=\textup{(\alph*)}]
\item the only edge entering a cyclic vertex is the preceding edge
of its cycle;
\item the noncyclic vertices induce an acyclic subgraph from which
no cycle is reachable;
\item every $(w,h)\in\mathcal H_d$ satisfies
$\min(|w|,|h|)\le a$.
\end{enumerate}
In particular, $h_{r,s}=0$ if $\min(r,s)>a$, and
\begin{equation}\label{gk:bounded-dimension}
 \dim_\kk\mathfrak g_d
 =\ell_{|d|}+\sum_{m=2}^{a}h_{m+|d|,m}<\infty,
\end{equation}
with an empty sum when $a<2$.
\end{lemma}
\begin{proof}
Every strongly connected component containing a cycle is a single
simple directed cycle. Indeed, otherwise some vertex has two
different outgoing edges within that component. Completing each edge
by a shortest return path produces two distinct simple cycles through
that vertex.

Contract the strongly connected components. The resulting graph is
a finite directed acyclic graph, and each of its vertices is reachable
from a source. A source component must contain a cycle: otherwise it
would be a single noncyclic vertex with no incoming edge, contrary to
the definition of $\Gamma_\circ$. If a cyclic component had an incoming
edge from another component, following predecessors back to a source
would give a path from another cyclic component to it. This is
forbidden. Likewise, after leaving a cycle a path cannot return to it
or reach another cycle. This proves \textup{(a)} and \textup{(b)}.

Let $u=s(w)=s(h)$. If $u$ is noncyclic, both paths stay in the
acyclic subgraph on $a$ vertices and have length at most $a-1$.
If $u$ is cyclic, the two paths cannot both start with the cycle edge:
that edge is the unique incoming edge of its range and hence is the
chosen edge $\lambda$ there. Such a common first edge is excluded in
the definition of $\mathcal H_d$. Thus at least one path leaves the
cycle immediately. Its remaining segment lies in the acyclic subgraph,
so its total length is at most $a$.
This proves \textup{(c)}. Formula~\eqref{gk:bounded-dimension} now
follows from Theorem~\ref{fin:dimension}, since a finite graph has
only finitely many paths of each fixed length.
\end{proof}

\subsection{Eventual periodicity and a uniform bound}
\label{gk:period-section}

\begin{proposition}\label{gk:periodicity}
Suppose condition~\textup{(iv)} holds and $\Gamma$ contains at least
one cycle. Let $\ell_1,\ldots,\ell_t$, with $t\ge1$, be its cycle
lengths. Put
\[
 p_\Gamma=\operatorname{lcm}(\ell_1,\ldots,\ell_t),
 \qquad d_0=\max(1,a),
\]
where $a$ is as in Lemma~\ref{gk:core-shape}. Then
\begin{equation}\label{gk:period}
 \dim_\kk\mathfrak g_{d+p_\Gamma}=\dim_\kk\mathfrak g_d
 \qquad(d\ge d_0).
\end{equation}
Moreover, $\dim\mathfrak g_{-d}=\dim\mathfrak g_d$, and these dimensions
have a common finite bound. The nonnegative part has the rational
generating series
\begin{equation}\label{gk:hilbert}
\begin{split}
 \sum_{d\ge0}(\dim_\kk\mathfrak g_d)t^d
 ={}&\sum_{d=0}^{d_0-1}(\dim_\kk\mathfrak g_d)t^d\\
 &+\frac{t^{d_0}}{1-t^{p_\Gamma}}
   \sum_{j=0}^{p_\Gamma-1}(\dim_\kk\mathfrak g_{d_0+j})t^j.
\end{split}
\end{equation}
If $\Gamma$ is acyclic, then $\mathfrak g=0$ and the series is zero.
\end{proposition}
\begin{proof}
Order the vertices of $\Gamma_\circ$ with cyclic vertices first.
Lemma~\ref{gk:core-shape} gives the block adjacency matrix
\[
 B=\begin{pmatrix}P&R\\0&N\end{pmatrix},
 \qquad P^{p_\Gamma}=I,\qquad N^a=0\quad(a>0),
\]
where $P$ is the direct sum of the permutation matrices of the cycles.
For $a=0$, simply $B=P$. If $a>0$ and $n\ge a$, then
\[
 B^n=\begin{pmatrix}
 P^n&\displaystyle\sum_{j=0}^{a-1}P^{n-1-j}RN^j\\
 0&0
 \end{pmatrix}.
\]
Consequently,
\begin{equation}\label{gk:matrix-period}
 B^{n+p_\Gamma}=B^n\qquad(n\ge a),
\end{equation}
also when $a=0$. This is an equality of integer matrices, independent
of the coefficient field.

For $d>0$, formula~\eqref{fin:low-rank} expresses $\ell_d$ in terms of
$B^{d+1},B^d$, and $\mu_d$. The latter number is also
$p_\Gamma$-periodic, since the length of every isolated cycle divides
$p_\Gamma$. For $2\le m\le a$, the first argument of
formula~\eqref{fin:mixed-count-formula} for $h_{m+d,m}$ involves the
powers $m+d,m+d-1,m+d-2$, all at least $a$ when $d\ge d_0$.
Thus every term of~\eqref{gk:bounded-dimension} is unchanged upon
replacing $d$ by $d+p_\Gamma$.

The symmetry of positive and negative degrees follows from
Theorem~\ref{fin:dimension}. The first $d_0+p_\Gamma$ dimensions are
finite, and all subsequent values repeat; their maximum is a uniform
bound. Summing the periodic tail gives~\eqref{gk:hilbert}.
\end{proof}

\begin{proof}[Proof of Theorem~\ref{gk:main}]
The implications \textup{(iii)}$\Rightarrow$\textup{(ii)}
$\Rightarrow$\textup{(i)} are immediate.
Lemma~\ref{gk:witnesses} proves
\textup{(i)}$\Rightarrow$\textup{(iv)}, while
Lemma~\ref{gk:core-shape} and Proposition~\ref{gk:periodicity} prove
\textup{(iv)}$\Rightarrow$\textup{(iii)}. For an acyclic graph,
$\mathfrak g=0$.

To relate the graph condition to growth, use
\cite[Corollary~3 and Theorem~5]{AAJZ}.
Intersecting distinct cycles imply exponential growth and infinite
GK-dimension. If distinct cycles are vertex-disjoint, then
\begin{equation}\label{gk:aajz}
 \operatorname{GKdim}A=\max\{2d_1-1,\,2d_2\},
\end{equation}
where $d_1$ is the maximum number of cycles in a directed chain and
$d_2$ is the maximum number in a chain whose final cycle has an exit.
An absent type of chain has length zero. In particular the formula
gives $\max(-1,0)=0$ for an acyclic graph.

Condition~\textup{(iv)} gives $d_1,d_2\le1$, hence
$\operatorname{GKdim}A\le2$. Conversely, if the GK-dimension is at
most two, cycles are vertex-disjoint and~\eqref{gk:aajz} forces
$d_1\le1$. Thus no path connects distinct cycles. This proves
\textup{(iv)}$\Leftrightarrow$\textup{(v)}.

When the conditions fail, Lemma~\ref{gk:witnesses} supplies a countably
infinite independent family in $\mathfrak g_0$. A finite graph has
at most countably many finite paths, so Theorem~\ref{fin:main}
supplies at most countably many basis elements. The dimension is
therefore exactly $\aleph_0$.
\end{proof}

\subsection{Effective verification and scope of the criterion}
\label{gk:scope}

Condition~\textup{(iv)} can be checked without listing all simple cycles.
Compute the strongly connected components and check that every cyclic
component has precisely one internal incoming and one internal
outgoing edge at each vertex. In the condensation graph, check that
no cyclic component is reachable from a different cyclic component.
For this last test, a topological traversal need only propagate the
Boolean information that a cyclic component has occurred upstream.
Thus the whole test takes linear time in the number of vertices and
edges. Computing the normal form of a particular class still uses
the procedure of Section~\ref{fin:algorithm-section}.

\begin{example}[The transition from GK-dimension two to three]
\label{gk:two-loops-example}
A loop without exits has GK-dimension one and
$\dim\mathfrak g_d=1$ for every $d$. A loop with one exit to a sink
has GK-dimension two and the same homogeneous dimensions.
Replacing the sink by a vertex with a second loop gives $d_1=2$,
$d_2=1$, hence GK-dimension three. If the loops are $a,b$ and the
connecting edge is $f$, then
\[
 [\cD_{s(a),(a^nf)(fb^n)^*}],\qquad n\ge1,
\]
are linearly independent classes of degree zero.
\end{example}

\begin{remark}[Scope]\label{gk:limits}
The criterion does not recover the precise GK-dimension from the
abstract Lie algebra $\mathfrak g$. The Laurent algebra and the
algebraic Toeplitz algebra have isomorphic graded outer-derivation
Lie algebras, as shown in Sections~\ref{app:cycles}
and~\ref{app:toeplitz}, but their GK-dimensions are one and two.

The use of outer classes is essential. For the Toeplitz graph with
loop $a$ and exit $b$, the degree-zero derivations
$\ad_{a^nbb^*(a^*)^n}$, $n\ge0$, are linearly independent.
Indeed, in a finite linear combination with largest nonzero index
$N$, its value on $a$ has a nonzero coefficient at the matrix unit
$a^{N+1}bb^*(a^*)^N$. Thus $\Der_\kk(A)_0$ is infinite-dimensional,
although $\dim\mathfrak g_0=1$.

Nor need every homogeneous component be infinite-dimensional when
the GK-dimension exceeds two. Join two cycles of length two by one
edge, choosing the bipartition so that the resulting graph is
bipartite. Its GK-dimension is three and
$\dim\mathfrak g_0=\aleph_0$. But any two parallel paths have lengths
of the same parity, so the corners $vAv$ have only even degrees and
$\mathfrak g_{2j+1}=0$.

Finiteness of the graph is also essential. For the incoming star in
Section~\ref{inf:short-paths}, the algebra is locally matricial and has
GK-dimension zero when the invariant for a non-finitely generated
algebra is defined as the supremum over its finitely generated
subalgebras. Nevertheless, its degree-zero outer derivations contain
the vector space
\[
 \kk^{\NN}/\bigl(\kk^{(\NN)}+\kk\mathbf1\bigr).
\]
This is infinite-dimensional: partition $\NN$ into countably many
infinite sets. The characteristic functions of these sets have
linearly independent classes, since a finite combination that is
constant outside a finite set must vanish on each part.
\end{remark}

Formula~\eqref{gk:aajz} is the established growth theorem of
\cite{AAJZ}. Its application to $\mathfrak g_0$ uses the independence
of the outer classes in Theorem~\ref{fin:main}, not only a count
of words in the algebra.

\section{Finite examples and explicit Lie brackets}
\label{app:examples}

The graph is finite throughout this section. Equalities of
derivations and descriptions of free modules hold over the
commutative coefficient ring $\kk$; dimension statements additionally
assume that $\kk$ is a field. For a closed real path $c$ based at $v$,
we retain the path-operator conventions
\[
 \cD_c=\cD_{v,c},\qquad
 \cD_{c^*}=\cD_c^\sharp=-\cD_{v,c^*},
 \qquad D^\sharp(x)=D(x^*)^*.
\]
Thus the original ghost-path operator has the opposite sign to the
operator with parameter $c^*$. Brackets around one derivation denote
its class in $\HH^1_\kk(A)$; brackets between two derivations denote
their ordinary commutator.

\subsection{Acyclic graphs and a source preceding a loop}
\label{app:acyclic}

Repeated deletion of sources exhausts a finite acyclic graph, so
Theorem~\ref{fin:main} gives $\HH^1_\kk(L_\kk(\Gamma))=0$.
Equivalently, its Leavitt path algebra is a finite direct sum of
full matrix algebras over $\kk$, whose $\kk$-linear derivations
are inner.

For a graph $u\xrightarrow{b}v$ with one loop $a$ at $v$, deleting
$u$ leaves the one-loop graph. Here
$L_\kk(\Gamma)\simeq M_2(\kk[t,t^{-1}])$.
Its outer derivations, with their degrees, agree with those of one
loop. Representatives are $\cD_{v,v}$, $\cD_{v,a^n}$, and
$\cD_{v,(a^*)^n}$ for $n\ge1$, with parameter zero at the deleted
vertex. The source-free counting formulas apply after deletion.

\subsection{One loop and an oriented cycle}
\label{app:cycles}

For one loop $a=t$, the algebra is $A=\kk[t,t^{-1}]$.
A derivation may take any value on $t$, and the identity $tt^{-1}=1$
then determines its value on $t^{-1}$. Hence
\begin{equation}\label{app:witt-basis}
 W_\kk:=\Der_\kk(A)
 =\bigoplus_{n\in\ZZ}\kk L_n,\qquad
 L_n=t^{n+1}\frac{d}{dt},\qquad
 [L_m,L_n]=(n-m)L_{m+n}.
\end{equation}
There are no nonzero inner derivations. We have
$L_n=\cD_{v,t^n}$ for all $n\in\ZZ$, including the Euler derivation
$L_0(t)=t$, $L_0(t^{-1})=-t^{-1}$.
For $n>0$, the original conventions read
$\cD_{a^n}=L_n$ and $\cD_{(a^n)^*}=-L_{-n}$.
In positive characteristic the integer coefficients in
\eqref{app:witt-basis} are interpreted in $\kk$; no simplicity
assertion is intended.

Now let the entire graph be a cycle of length $r$, with base vertex
$v_0$ and based cycle $c$. Then $v_0Av_0=\kk[c,c^{-1}]$.
For each vertex $v_i$, let $q_i$ be the unique path from $v_i$ to
$v_0$ of length less than $r$, with $q_0=v_0$. The relations give
\[
 q_i^*q_j=\delta_{ij}v_0,\qquad q_iq_i^*=v_i.
\]
The map $E_{ij}t^n\mapsto q_ic^nq_j^*$ identifies
$A$ with $M_r(\kk[t,t^{-1}])$.

For completeness, if $D$ is a derivation of $M_r(R)$ and
$h=\sum_iD(E_{i1})E_{1i}$, then
$[h,E_{kl}]=D(E_{kl})$.
The derivation $D-\ad_h$ fixes all matrix units and acts on every
entry by the same derivation of $R$. When $R$ is commutative, an
inner derivation vanishes on the centre $RI_r$. Therefore
$\HH^1_\kk(M_r(R))\simeq\Der_\kk(R)$ as Lie algebras, without
division by $r$.

In graph coordinates, choose
\[
 D_n=\cD_{v_0,c^n},\qquad n\in\ZZ.
\]
Negative powers are taken in the corner with identity $v_0$.
If $e$ is the first edge of $c$, then
$D_n(e)=c^ne$, $D_n(e^*)=-e^*c^n$; all other edges and ghost edges
are fixed. Thus $D_n(c)=c^{n+1}$, and
\begin{equation}\label{app:cycle-witt}
 [D_m,D_n]=(n-m)D_{m+n},\qquad \deg D_n=nr.
\end{equation}
This is an equality in $\Der_\kk(A)$. The paths $q_i$ are fixed by
every $D_n$, so the matrix identification above shows that the
classes $[D_n]$ give a basis of $\HH^1_\kk(A)$.
Over a field,
\[
 \dim_\kk\HH^1_\kk(A)_d=
 \begin{cases}1,&r\mid d,\\0,&r\nmid d.\end{cases}
\]
For the two-cycle $a:u\to v$, $b:v\to u$, one also has
\[
 \cD_{u,aa^*}-\cD_{v,bb^*}=\ad_u.
\]
Its two vertices therefore do not give two independent outer
directions of degree zero.

\subsection{A loop with one exit: the algebraic Toeplitz algebra}
\label{app:toeplitz}

Let $a$ be a loop at $v$ and $b:v\to w$ an exit to a sink.
Put $J=AwA$. A basis of $vAv$ is
\[
 v;\qquad a^n,(a^*)^n\ (n\ge1);\qquad
 E_{ij}=a^ibb^*(a^*)^j\ (i,j\ge0).
\]
Here $a^*a=v$, $aa^*=v-bb^*$, and
$E_{ij}E_{kl}=\delta_{jk}E_{il}$.
Thus $I=vJv=\bigoplus_{i,j\ge0}\kk E_{ij}$ and
$vAv/I\simeq\kk[t,t^{-1}]$.

The vertex $v$ is the only regular vertex, so a normalized derivation
is determined by $p\in vAv$.
For a diagonal implementer $h+\lambda w$, its parameter is
\begin{equation}\label{app:toeplitz-boundary}
 h-aha^*-\lambda bb^*.
\end{equation}
Write $T(h)=h-aha^*$. Direct calculation gives
\[
 T(v)=E_{00},\quad T(a^n)=E_{n0},\quad
 T((a^*)^n)=E_{0n},\quad T(E_{ij})=E_{ij}-E_{i+1,j+1}.
\]
Hence $\im T=I$: the inclusion $\im T\subseteq I$ follows after
passing to the Laurent quotient, and the reverse inclusion follows
by finite induction on $\min(i,j)$. The term $\lambda bb^*$ is
already in $\im T$.

Every normalized derivation preserves $J$, since it vanishes on
$w$, and induces $t\mapsto\overline p\,t$ on
$A/J\simeq\kk[t,t^{-1}]$. The kernel of this parameter map is $I$,
precisely the space of inner parameters.
Every derivation differs from a normalized one by an inner
derivation. We therefore obtain an isomorphism of Lie algebras
\begin{equation}\label{app:toeplitz-witt}
 \HH^1_\kk(A)\simeq W_\kk.
\end{equation}
Each homogeneous component is free of rank one.

To express the calculation in the original operators, put
\[
 H=\cD_{aa^*},\qquad P_n=\cD_{a^n},\qquad
 G_n=\cD_{(a^n)^*}=P_n^\sharp\qquad(n\ge1).
\]
They vanish on vertices and have the following values:
\[
\begin{array}{c|cccc}
 &a&b&a^*&b^*\\ \hline
 H&a&0&-a^*&0\\
 P_n&a^{n+1}&a^nb&-a^{n-1}&0\\
 G_n&-(a^*)^{n-1}&0&(a^*)^{n+1}&b^*(a^*)^n
\end{array}
\qquad(a^0=(a^*)^0=v).
\]
The Leibniz rule gives the exact identities
\begin{align}
 [H,P_n]&=nP_n,&[H,G_n]&=-nG_n,\notag\\
 [P_n,P_m]&=(m-n)P_{n+m},&
 [G_n,G_m]&=(m-n)G_{n+m}.
 \label{app:toeplitz-pure-brackets}
\end{align}
Under~\eqref{app:toeplitz-witt}, $[H],[P_n],[G_n]$ correspond to
$L_0,L_n,-L_{-n}$, respectively. Thus the mixed bracket in cohomology is
\begin{equation}\label{app:toeplitz-mixed-bracket}
 [[P_n],[G_m]]=
 \begin{cases}
 (n+m)[P_{n-m}],&n>m,\\
 2n[H],&n=m,\\
 -(n+m)[G_{m-n}],&n<m.
 \end{cases}
\end{equation}
In the full derivation algebra there can be an inner correction.
For example,
\begin{equation}\label{app:toeplitz-inner-correction}
 [P_1,G_1]=2H+\ad_v:
\end{equation}
both sides take the values $2a,b,-2a^*,-b^*$ on
$a,b,a^*,b^*$, respectively.
Other useful inner identities are
\begin{align*}
 \cD_{bb^*}&=\ad_v=-\ad_w,\\
 \cD_{v,a^n}-\cD_{v,a^{n+1}a^*}&=\ad_{a^n},\\
 \cD_{v,(a^*)^n}-\cD_{v,a(a^*)^{n+1}}&=\ad_{(a^*)^n}
 \qquad(n\ge1).
\end{align*}
Consequently, the corresponding pure and mixed parameters do not
give two independent families. The Euler derivation
$\cD_{v,v}=H+\cD_{bb^*}$ has the same outer class as $H$.
These identities correct the free-family assertion in
\cite[Section~6, p.~82]{Lopatkin}.
There is also a sign correction in the positive mixed family of
\cite[Theorem~6.1, p.~84]{Lopatkin}. If
$Q_n=\cD_{v,a^na^*}$ for $n\ge1$, then
\[
Q_n(a)=a^n,\qquad
[Q_n,Q_m]=(m-n)Q_{n+m-1}.
\]
This follows by evaluation on $a,b,a^*,b^*$, with
$Q_n(a^*)=-a^{n-1}a^*$ and $Q_n(b)=Q_n(b^*)=0$.

\subsection{One loop with parallel exits}
\label{app:parallel-exits}

Let $a$ be a loop at $v$, and let $b_{ij}:v\to w_i$ be exits to
distinct sinks $w_1,\ldots,w_q$, where $1\le j\le m_i$ and $m_i\ge1$.
Put $S=\sum_i m_i^2$. A basis of $vAv$ consists of
$v,a^n,(a^*)^n$ together with
\[
 X_{i;jk}^{nm}=a^nb_{ij}b_{ik}^*(a^*)^m,
 \qquad n,m\ge0,\quad 1\le j,k\le m_i.
\]
The full diagonal inner relation is
\[
 h-aha^*-\sum_i\lambda_i\sum_jb_{ij}b_{ij}^*.
\]
Again let $T(h)=h-aha^*$. In degree $d$ the relations
\[
 T(X_{i;jk}^{nm})=X_{i;jk}^{nm}-X_{i;jk}^{n+1,m+1}
\]
identify all symbols with the same $i,j,k$ and $n-m=d$.
Their quotient has one copy of $M_{m_i}(\kk)$ for each $i$, and
one pure parameter of degree $d$.
For $d>0$, the remaining relation $T(a^d)$ is the sum of the
identity matrices of all blocks; for $d<0$, the same is true of
$T((a^*)^{-d})$. In degree zero, the sink parameters give the
identity matrix of each block separately, and $T(v)$ is their sum.
These are all the relations, because they are the images under
$\partial$ of the stated basis of the diagonal corners.
Each elimination has a unit coefficient. Hence, over a field,
\begin{equation}\label{app:exits-dimensions}
 \dim_\kk\HH^1_\kk(A)_d=
 \begin{cases}
 S-q+1,&d=0,\\
 S,&d\ne0.
 \end{cases}
\end{equation}
Over a commutative ring the same argument gives free modules of
these ranks.

The derivations $T_{i;jk}=\cD_{v,b_{ij}b_{ik}^*}$ act as matrix
units on the parallel real edges and as their negative transposes
on the ghost edges. Thus
\[
 [T_{i;jk},T_{i;lm}]
 =\delta_{kl}T_{i;jm}-\delta_{mj}T_{i;lk},
\]
and operators from different blocks commute.
The identity of block $i$ is the inner derivation $-\ad_{w_i}$.
The Euler derivation $\cD_{v,v}$ commutes with all degree-zero
derivations. The preceding complete list of relations therefore gives
\begin{equation}\label{app:exits-zero-lie}
 \HH^1_\kk(A)_0\simeq
 \kk\oplus\bigoplus_{i=1}^{q}
 \bigl(\gl_{m_i}(\kk)/\kk I_{m_i}\bigr)
\end{equation}
as Lie algebras, with the first summand central.
If $\kk$ is a field of characteristic zero, the matrix quotients
are isomorphic to $\slie_{m_i}(\kk)$.
The quotient formulation is retained in arbitrary characteristic.

Two exits to the same sink give dimension four in every degree.
One exit to each of two different sinks gives dimension one in
degree zero and two in every nonzero degree. Thus the grouping of
edge ranges matters. Formulas~\eqref{app:exits-dimensions}
and~\eqref{app:exits-zero-lie} specify the homogeneous dimensions
and the degree-zero Lie algebra; the remaining brackets are computed
using~\eqref{gen:normalizedbracket} and Theorem~\ref{fin:main}.

\subsection{Bouquets: matrix derivations and path insertions}
\label{app:bouquet}

Consider one vertex with $q\ge2$ loops $e_1,\ldots,e_q$.
The operators $T_{ij}=\cD_{e_ie_j^*}$ satisfy
\begin{equation}\label{app:bouquet-matrices}
 [T_{ij},T_{kl}]
 =\delta_{jk}T_{il}-\delta_{li}T_{kj}.
\end{equation}
Their outer classes are independent: the remaining low-depth
relations of degree zero come from vertices, and the relation for
the single vertex is zero.
Thus $\HH^1_\kk(A)_0$ contains $\gl_q(\kk)$.
The sum $\mathcal E=\sum_iT_{ii}=\cD_{v,v}$ is the nonzero Euler
derivation.

For positive words $c=c_1\cdots c_n$ and $d=d_1\cdots d_m$, direct
application of the Leibniz rule gives
\begin{equation}\label{app:bouquet-insertions}
\begin{split}
 [\cD_c,\cD_d]
 ={}&\sum_{j=2}^{m}
 \cD_{d_1\cdots d_{j-1}c\,d_j\cdots d_m}\\
 &-\sum_{i=2}^{n}
 \cD_{c_1\cdots c_{i-1}d\,c_i\cdots c_n}.
\end{split}
\end{equation}
Indeed, $\cD_c(d)$ inserts $c$ before each letter of $d$, and
$\cD_d(c)$ inserts $d$ before each letter of $c$.
The first insertions cancel the term $dc-cd$ in the parameter bracket.
Every word on the right has length $n+m$; in particular,
$[\cD_a,\cD_b]=0$ for any two loops $a,b$.
For the original ghost operator,
\begin{equation}\label{app:bouquet-one-letter-mixed}
 [\cD_a,\cD_b^\sharp]=\cD_{ab^*}+\delta_{ab}\mathcal E.
\end{equation}
This follows from~\eqref{gen:normalizedbracket} with parameters
$a,-b^*$ and the identity $b^*a=\delta_{ab}v$.

For $n>0$, there are $q^{n+2}$ low-depth pairs $(P,e)$ with
$|P|=n+1$ and $q^n$ independent relations from closed paths of
length $n$. Hence
\[
 \ell_n=(q^2-1)q^n\quad(n>0),\qquad \ell_0=q^2.
\]
For lengths $r,s\ge2$, the restrictions on the first and last edge
pairs are independent, and
\begin{equation}\label{app:bouquet-deep-count}
 h_{r,s}=(q^2-1)^2q^{r+s-4}.
\end{equation}
Theorem~\ref{fin:main} gives independent outer classes for all these
pairs. Keeping $r-s=d$ fixed and letting their lengths increase yields
\[
 \dim_\kk\HH^1_\kk(L_\kk(1,q))_d=\aleph_0
 \qquad(d\in\ZZ)
\]
over any field.

For two loops $a,b$, choose $a$ as both the special outgoing and the
selected incoming edge. In degree one, order the columns as
\[
 (aa,a),(ab,a),(ba,a),(bb,a),(aa,b),(ab,b),(ba,b),(bb,b).
\]
The two inner relations form the matrix
\begin{equation}\label{app:bouquet-matrix-one}
 M_1=\begin{pmatrix}
 0&0&0&0&0&1&-1&0\\
 0&-1&1&0&0&0&0&0
 \end{pmatrix}.
\end{equation}
The rows correspond to $a,b$, namely
$(ab,b)-(ba,b)$ and $(ba,a)-(ab,a)$.
One may eliminate columns $(ab,a),(ab,b)$, keeping the other six,
in addition to the infinite mixed family.
In degree zero,
\[
 [\cD_{ab^*},\cD_{ba^*}]
 =\cD_{aa^*}-\cD_{bb^*}.
\]
To distinguish normal words in $A$ from representatives in
cohomology, observe that
\[
 \cD_{bb^*}-\cD_{ab(ab)^*}-\cD_{bb(bb)^*}=\ad_{bb^*}.
\]
Thus
$[\cD_{ab(ab)^*}]=[\cD_{bb^*}]-[\cD_{bb(bb)^*}]$.
The word on the left is normal in $A$, since its last edges are the
nonspecial edge $b$, but its common first edge $a$ is eliminated
in the outer normal form.

\subsection{Two loops joined by an edge}
\label{app:linked-cycles}

Let $a$ be a loop at $v$, let $c$ be a loop at $w$, and let
$b:v\to w$. Then
\[
 B=\begin{pmatrix}1&1\\0&1\end{pmatrix},
 \qquad
 B^n=\begin{pmatrix}1&n\\0&1\end{pmatrix}.
\]
For $n>0$ there are $n+3$ low-depth coordinates: one each for
$a,c$, and $n+1$ for $b$.
The independent relations from $a^n,c^n$ eliminate
$(a^nb,b)$ and $(bc^n,b)$, respectively. Consequently,
\[
 \ell_n=n+1\quad(n>0),\qquad \ell_0=2.
\]
In degree zero, the three coordinates $(a,a),(b,b),(c,c)$ have
one vertex relation; the positive-degree formula therefore does
not extend to $n=0$.
For $n>0$ the low-depth part can be represented by
\[
 [\cD_{a^n}],\qquad [\cD_{c^n}],\qquad
 [\cD_{v,a^ibc^{n-i}b^*}]\quad(1\le i\le n-1).
\]
The final list is empty for $n=1$, and the negative-degree part is
obtained by applying $\sharp$.

Choose $a,c$ as the special outgoing and selected incoming edges.
For each $r,s\ge2$ there are exactly two remaining mixed pairs:
\begin{equation}\label{app:linked-deep-pairs}
 (a^{r-1}b,\ bc^{s-1}),\qquad
 (bc^{r-1},\ a^{s-1}b).
\end{equation}
Indeed, to avoid a common special final edge $c$, one path must end
with $b$; to avoid a common selected initial edge $a$, the other
must start with $b$. These choices give exactly the displayed pairs.
Pairs ending at $v$ have the forbidden common final edge $a$.
By Theorem~\ref{fin:main}, the operators with these parameters give
independent outer classes. There are infinitely many at every fixed
difference $r-s$. Every homogeneous component is therefore
countably infinite-dimensional over a field.

\subsection{Homogeneous dimensions}
\label{app:table}

The table records dimensions over a field. Over a nonzero commutative
ring, the same cardinal numbers are the ranks of the free modules
constructed above.
\begin{center}
\begin{tabular}{p{0.43\textwidth}p{0.17\textwidth}p{0.27\textwidth}}
\toprule
Graph & $\dim\HH^1_0$ & $\dim\HH^1_d$, $d\ne0$\\
\midrule
Finite acyclic graph & $0$ & $0$\\
One loop & $1$ & $1$\\
Cycle of length $r$ & $1$ & $1$ if $r\mid d$, otherwise $0$\\
Loop with one exit to a sink & $1$ & $1$\\
Loop with $m_i$ exits to sink $i$ & $\sum_i m_i^2-q+1$ & $\sum_i m_i^2$\\
Bouquet of $q\ge2$ loops & $\aleph_0$ & $\aleph_0$\\
Two loops joined by an edge & $\aleph_0$ & $\aleph_0$\\
\bottomrule
\end{tabular}
\end{center}
The full Lie algebra is isomorphic to $W_\kk$ for a cycle and for
a loop with one exit; for the cycle its grading is rescaled by the
cycle length. In the other examples, the table records dimensions,
and the brackets are determined by the explicit path operators and
reduction of inner relations.

\section{Infinite graphs with short paths}
\label{inf:short-paths}

\subsection{The incoming star and finitary matrices}
Let $V=\{v_0,v_1,\ldots\}$, with $v_0$ a sink, and let
$e_n:v_n\to v_0$ for $n\ge1$. Every nontrivial path has length one.
Set $p_0=v_0$ and $p_n=e_n$ for $n\ge1$.
The left and right supports of a ghost edge are
\begin{equation}\label{inf:star-supports}
v_0e_n^*=e_n^*=e_n^*v_n,\qquad
e_n^*v_0=e_n^*v_nv_0=0.
\end{equation}
Also $v_0e_n=v_0v_ne_n=0$. Consequently
$p_i^*p_j=\delta_{ij}v_0$ for all $i,j\ge0$, including the cases
in which an index is zero. The elements $F_{ij}=p_ip_j^*$ satisfy
\[
F_{ij}F_{kl}=\delta_{jk}F_{il}.
\]

For completeness, the isomorphism with the finitary matrix algebra
is explicit in both directions. The formulas
\begin{equation}\label{inf:star-map}
\rho(v_i)=E_{ii},\qquad
\rho(e_n)=E_{n0},\qquad
\rho(e_n^*)=E_{0n}
\end{equation}
define a homomorphism $\rho:A\to M_\fin(\mathbb N_0,\kk)$.
The support relations and orthogonality are immediate, and
\[
E_{0n}E_{m0}=\delta_{nm}E_{00},\qquad
E_{n0}E_{0n}=E_{nn}
\]
verify the remaining defining relations. Conversely,
$\psi(E_{ij})=F_{ij}$ extends to a homomorphism on finite linear
combinations of matrix units. We have $\rho\psi=\mathrm{id}$ and
\[
F_{00}=v_0,\qquad F_{n0}=e_n,\qquad F_{0n}=e_n^*,\qquad F_{nn}=v_n,
\]
so $\psi\rho$ is the identity on the generators of $A$. Thus
\[
A\simeq M_\fin(\mathbb N_0,\kk),\qquad
\cM\simeq M_\rcf(\mathbb N_0,\kk).
\]
Here $M_\fin$ means finite total support, while $M_\rcf$ means
finite support in every row and every column. The same construction
for a finite star with $N$ edges gives $M_{N+1}(\kk)$.

\begin{proposition}\label{inf:star}
For the incoming star there are isomorphisms of $\kk$-Lie algebras
\begin{equation}\label{inf:star-quotients}
\Der_\kk(A)\simeq
\frac{M_\rcf(\mathbb N_0,\kk)}{\kk I},
\qquad
\HH^1_\kk(A)\simeq
\frac{M_\rcf(\mathbb N_0,\kk)}
{M_\fin(\mathbb N_0,\kk)+\kk I}.
\end{equation}
\end{proposition}
\begin{proof}
Each corner $v_nAv_n=\kk v_n$ has rank one. For the diagonal family
$(a_nv_n)_{n\ge0}$, the $n$th component of $\partial$ is
$(a_n-a_0)v_n$ for $n\ge1$. Hence $\partial:\cQ\to\cP$ is
surjective: take $a_0=0$ and choose the remaining $a_n$ freely.
The full-coordinate theorem now shows that every derivation is
a commutator with a matrix in $M_\rcf$. A matrix commuting with
all matrix units is scalar, giving the first kernel $\kk I$.
Taking the quotient by ordinary inner derivations gives the second
formula.
\end{proof}

This is a special case of the known description of derivations of
infinite matrix rings; see \cite[Theorem~1(a)]{Bez}. The argument
above obtains it directly from graph coordinates. Over $\mathbb C$,
the normalized derivation
\[
D(v_n)=0,\qquad D(e_n)=n e_n,\qquad D(e_n^*)=-n e_n^*
\]
is outer. Its implementing diagonal matrix is not the sum of a
finitary matrix and a scalar matrix. Thus bounded path length does
not imply $\HH^1=0$ for an infinite row-finite graph.

\subsection{An outer elementary operator with bounded path length}
Consider edges $e,f:v\to w$, where $w$ is a sink, together with
$g_n:u_n\to v$, where the $u_n$ are distinct sources.
Every path has length at most two. The elementary operator
\[
D=\cD_{ef^*},\qquad D(f)=e,\qquad D(e^*)=-f^*,
\]
vanishes on all other generators and is outer.

Indeed, if $D=\ad_h$, then $h$ commutes with all vertices, so
$h=\sum_{x\in F}h_x$, where $F$ is finite and $h_x\in xAx$.
Choose $u_n\notin F$. Then $0=D(g_n)=-g_nh_v$; multiplying on the
left by $g_n^*$ gives $h_v=0$. Since $wAw=\kk w$, the value
$[h,f]=-fh_w$ is a scalar multiple of $f$, contradicting $D(f)=e$
and the linear independence of real paths.

Consequently, even the innerness of a single $\cD_{wh^*}$ cannot
be characterized merely by the existence of arbitrarily long
incoming paths. The general criterion uses $\partial$ together
with finite support of the implementing element.

\section{Parallel layers and products of matrix Lie algebras}
\label{inf:layers}

Fix $r\ge2$. Let $V=\{v_n:n\ge0\}$, with exactly $r$ edges
$e_{n,1},\ldots,e_{n,r}$ from $v_n$ to $v_{n+1}$ for each $n$.
Put
\[
X_{ij}^{(n)}=\cD_{e_{n,i}e_{n,j}^*}.
\]
Their values on generators are
\[
X_{ij}^{(n)}(e_{m,k})=\delta_{nm}\delta_{jk}e_{n,i},\qquad
X_{ij}^{(n)}(e_{m,k}^*)=-\delta_{nm}\delta_{ik}e_{n,j}^*.
\]
It follows that
\begin{equation}\label{inf:layers-bracket}
[X_{ij}^{(n)},X_{kl}^{(m)}]
=\delta_{nm}\bigl(\delta_{jk}X_{il}^{(n)}-
\delta_{li}X_{kj}^{(n)}\bigr).
\end{equation}
In particular, different layers commute.

\begin{theorem}\label{inf:layers-theorem}
There are injective homomorphisms of $\kk$-Lie algebras
\begin{equation}\label{inf:layers-embeddings}
\prod_{n\ge0}\gl_r(\kk)\hookrightarrow\Der_\kk(A),
\qquad
\frac{\prod_{n\ge0}\gl_r(\kk)}
{\bigoplus_{n\ge0}\gl_r(\kk)}
\hookrightarrow\HH^1_\kk(A).
\end{equation}
The second map describes a Lie subalgebra of $\HH^1_\kk(A)$;
surjectivity is not asserted.
\end{theorem}
\begin{proof}
A family $M=(M_n)$ acts by
\[
D_M(e_{n,k})=\sum_i(M_n)_{ik}e_{n,i},\qquad
D_M(e_{n,k}^*)=-\sum_j(M_n)_{kj}e_{n,j}^*,\qquad
D_M(v_n)=0.
\]
These are finite sums for every generator, and
Theorem~\ref{gen:free} gives a derivation. Linear independence of
the edges proves injectivity, while \eqref{inf:layers-bracket}
proves that brackets are preserved.

Suppose first that only $M_n=M$ is nonzero. Define
\[
z_n=\sum_{i,j}M_{ij}e_{n,i}e_{n,j}^*,\qquad
h^{(n)}=\sum_{k=0}^{n}\ \sum_{\alpha:v_k\to v_n}
\alpha z_n\alpha^*.
\]
The trivial path is included for $k=n$. For fixed $k\le n$
there are exactly $r^{n-k}$ paths from $v_k$ to $v_n$.
Thus $h^{(n)}$ is a finite sum in $A$, supported at
$v_0,\ldots,v_n$. Write $h_k=v_kh^{(n)}v_k$. For $k<n$,
\[
h_k=\sum_i e_{k,i}h_{k+1}e_{k,i}^*.
\]
Consequently $[h^{(n)},e_{k,i}]=0$ before layer $n$, whereas
\[
[h^{(n)},e_{n,i}]=\sum_jM_{ji}e_{n,j}.
\]
The values after layer $n$ vanish. The ghost-edge values agree
with the displayed formula for $D_M$, so $D_M=\ad_{h^{(n)}}$.
Every finitely supported family therefore gives an ordinary inner
derivation.

Conversely, suppose $D_M=\ad_h$ with $h\in A$.
Since $D_M$ vanishes on vertices, $h$ is diagonal with respect to
the vertices and has finite vertex support. For all sufficiently
large $n$, both $v_nh$ and $v_{n+1}h$ vanish, so
$[h,e_{n,i}]=0$ for every $i$. Linear independence of the edges
implies $M_n=0$. Hence the intersection of the image of the product
with $\Inn(A)$ is exactly the image of the direct sum, proving
the second injection.
\end{proof}

This example exhibits the difference between admissible pointwise
finite families of path operators and their finite linear combinations.
\begin{proposition}[Multiplier implementation of layer derivations]
\label{inf:layers-multipliers}
For the parallel-layer graph of Theorem~\ref{inf:layers-theorem}, let
$D_M$ correspond to a family $M=(M_n)_{n\ge0}$ in
$\prod_{n\ge0}\gl_r(\kk)$. Then
\begin{equation}\label{inf:layers-multiplier-criterion}
D_M=\ad_b|_A\text{ for some }b\in\cM
\quad\Longleftrightarrow\quad
M_n\in\kk I_r\text{ for all sufficiently large }n.
\end{equation}
Every $D_M$ is locally inner: on each prescribed finite subset of
$A$ it agrees with an ordinary inner derivation. In particular, the
constant family $M_n=E_{12}$ gives a locally inner derivation that
cannot be implemented by a multiplier.
\end{proposition}
\begin{proof}
For $N\ge0$, let
\[
B_N=\operatorname{span}_{\kk}
\{\alpha\beta^*: \alpha,\beta:v_0\longrightarrow v_N\}
\subseteq v_0Av_0.
\]
The paths from $v_0$ to $v_N$ are indexed by their $N$ edge labels,
so there are $r^N$ of them. Their products $\alpha\beta^*$ are
matrix units and are linearly independent: multiplying a relation
on the left by a fixed $\alpha^*$ and on the right by a fixed
$\beta$ recovers its corresponding scalar coefficient times
$v_N$. Thus $B_N\simeq M_{r^N}(\kk)$, with $B_0=\kk v_0$.
The Cuntz--Krieger relation at $v_N$ gives
\[
\alpha\beta^*=\sum_{i=1}^r
(\alpha e_{N,i})(\beta e_{N,i})^*.
\]
Under the path-label identification, the inclusion $B_N\subseteq
B_{N+1}$ is therefore $X\mapsto X\otimes I_r$.
Every monomial in $v_0Av_0$ has the form $\alpha\beta^*$ with
paths of the same length: their terminal vertices agree, and a
path from $v_0$ reaches $v_N$ exactly after $N$ edges. Refining
finitely many such monomials to a common length proves
\[
v_0Av_0=\bigcup_{N\ge0}B_N.
\]

For $L\ge1$, the restriction of $D_M$ to $B_L$ is the commutator
with the matrix
\[
K_L=\sum_{j=0}^{L-1}
I_r^{\otimes j}\otimes M_j\otimes
I_r^{\otimes(L-1-j)}.
\]
Indeed, applying the Leibniz rule to a path changes its $j$th
edge label according to $M_j$, while application to the ghost
path gives the corresponding negative right action. In particular,
for $X\in M_r(\kk)$ define
\[
T_n(X)=I_r^{\otimes n}\otimes X\in B_{n+1}.
\]
The earlier tensor factors commute with $T_n(X)$, so
\begin{equation}\label{inf:layers-tail-action}
D_M(T_n(X))=T_n([M_n,X]).
\end{equation}

Suppose $D_M=\ad_b|_A$ for $b\in\cM$. Since $D_M$ vanishes on
all vertices, $b$ commutes with them and is therefore diagonal
in the corner-matrix decomposition. Its component
$b_0=v_0bv_0$ lies in $v_0Av_0$, hence in some $B_N$.
For $n\ge N$, its image in $B_{n+1}$ acts on the first $N$
tensor factors and commutes with $T_n(X)$. It follows from
\eqref{inf:layers-tail-action} that
$T_n([M_n,X])=0$ for all $X\in M_r(\kk)$. The map $T_n$ is
injective, so $M_n$ commutes with every matrix. Over a commutative
ring the centralizer of the matrix units in $M_r(\kk)$ is
$\kk I_r$: commuting with the diagonal matrix units eliminates
the off-diagonal entries, and commuting with $E_{ij}$ identifies
the diagonal entries. Hence $M_n\in\kk I_r$ for $n\ge N$.

Conversely, suppose $M_n$ is scalar for all sufficiently large
$n$. Write $M=F+(\lambda_n I_r)_{n\ge0}$ with $F$ finitely
supported. By Theorem~\ref{inf:layers-theorem}, $D_F=\ad_h$ for
some $h\in A$. Define scalars recursively by
\[
c_0=0,\qquad c_{n+1}=c_n-\lambda_n,
\]
and let $c=\operatorname{diag}(c_nv_n)\in\cQ\subseteq\cM$.
Then
\[
[c,e_{n,i}]=(c_n-c_{n+1})e_{n,i}
=\lambda_ne_{n,i},\qquad
[c,e_{n,i}^*]=-\lambda_ne_{n,i}^*.
\]
Both derivations vanish on the vertices, so
$D_M=\ad_{h+c}|_A$. This proves
\eqref{inf:layers-multiplier-criterion} without dividing by $r$.

Finally, each finite subset of $A$ can be expressed using only
finitely many edge layers. Truncate $M$ beyond all these layers.
The resulting derivation agrees with $D_M$ on the chosen subset
and is inner by Theorem~\ref{inf:layers-theorem}.
The family $M_n=E_{12}$ is nowhere scalar for $r\ge2$, so the
criterion shows that its locally inner derivation is not a
multiplier commutator.
\end{proof}

\Needspace{20\baselineskip}
\begin{corollary}\label{inf:layers-exact-sequence}
Identify
\[
S=\frac{\prod_{n\ge0}\gl_r(\kk)}
{\bigoplus_{n\ge0}\gl_r(\kk)}
\]
with its image in $\HH^1_\kk(A)$ from
Theorem~\ref{inf:layers-theorem}. For the map $\jmath$ in
\eqref{inf:exact-equation}, scalar families give an isomorphism
of abelian $\kk$-Lie algebras
\begin{equation}\label{inf:layers-intersection}
S\cap\im\jmath\simeq
\frac{\prod_{n\ge0}\kk}{\bigoplus_{n\ge0}\kk}.
\end{equation}
Writing $\mathfrak l_r=\gl_r(\kk)/\kk I_r$, the restriction of
$\pi$ induces an injective Lie homomorphism
\begin{equation}\label{inf:layers-right}
\frac{\prod_{n\ge0}\mathfrak l_r}
{\bigoplus_{n\ge0}\mathfrak l_r}
\hookrightarrow\frac{\cP}{\partial\cQ}.
\end{equation}
\end{corollary}
\begin{proof}
Since every element of $A$ is a multiplier, $[D_M]\in\im\jmath$
if and only if $D_M$ itself is a multiplier commutator. By
Proposition~\ref{inf:layers-multipliers}, its family is then
scalar outside a finite set. Modulo finitely supported families
this is represented by a scalar family, and such a family has
zero class in $S$ exactly when its scalar sequence is finitely
supported. This proves \eqref{inf:layers-intersection}.

Consider the Lie homomorphism
\[
\prod_{n\ge0}\gl_r(\kk)\longrightarrow\cP/\partial\cQ,
\qquad M\longmapsto\pi([D_M]).
\]
Its kernel, by exactness and the criterion, consists exactly of
the families that are scalar outside a finite set. This is also
the kernel of the natural surjection
\[
\prod_{n\ge0}\gl_r(\kk)\longrightarrow
\frac{\prod_{n\ge0}\mathfrak l_r}
{\bigoplus_{n\ge0}\mathfrak l_r}.
\]
For completeness, representatives of classes in $\mathfrak l_r$
can be chosen with $(1,1)$ entry zero, by replacing $X$ with
$X-X_{11}I_r$. Thus the displayed map is surjective over every
commutative base ring, with no invertibility condition on $r$.
Taking the common quotient proves \eqref{inf:layers-right}.
\end{proof}

The layer families therefore exhibit nonzero classes in both
extreme terms of \eqref{inf:exact-equation}. They also separate
local innerness from implementation by a multiplier: every
$D_M$ is locally inner, while a nonscalar tail can contribute
nontrivially to $\cP/\partial\cQ$.

\appendix

\section{Infinitesimal, formal, and path automorphisms}
\label{aut:section}

In this section $A=L_\kk(\Gamma)$ for a finite nonempty graph.
Sources are allowed except where explicitly excluded. Automorphisms
are $\kk$-algebra automorphisms; they need not preserve the
anti-involution $*$. The identity is $1=\sum_{v\in V}v$.

\subsection{First-order automorphisms and outer derivations}
\label{aut:dual-numbers}

Let $R=\kk[\varepsilon]/(\varepsilon^2)$. Write $\mathcal G_1(A)$
for the group of $R$-automorphisms of $A\otimes_\kk R$ that are
the identity modulo $\varepsilon$, and $\mathcal I_1(A)$ for the
subgroup of conjugations by $1+\varepsilon h$, $h\in A$.

\begin{proposition}\label{aut:tangent}
The assignment $D\mapsto\mathrm{id}+\varepsilon D$ gives canonical
isomorphisms of additive groups
\[
\mathcal G_1(A)\simeq\Der_\kk(A),\qquad
\mathcal G_1(A)/\mathcal I_1(A)\simeq\HH^1_\kk(A).
\]
\end{proposition}
\begin{proof}
The decomposition $A\otimes_\kk R=A\oplus\varepsilon A$ writes
every such $R$-linear map uniquely as $\mathrm{id}+\varepsilon D$.
Multiplicativity is equivalent to the Leibniz rule, and the inverse
is $\mathrm{id}-\varepsilon D$. Composition corresponds to addition.
Finally,
\begin{equation}\label{aut:inner-first-order}
(1+\varepsilon h)x(1-\varepsilon h)=x+\varepsilon[h,x],
\end{equation}
so the subgroup $\mathcal I_1(A)$ corresponds exactly to $\ad(A)$.
\end{proof}

Over $\kk[\varepsilon,\eta]/(\varepsilon^2,\eta^2)$, direct
multiplication yields
\begin{equation}\label{aut:group-commutator}
\begin{split}
&(\mathrm{id}+\varepsilon D)(\mathrm{id}+\eta F)
(\mathrm{id}-\varepsilon D)(\mathrm{id}-\eta F)\\
&\hspace{25mm}=\mathrm{id}+\varepsilon\eta[D,F].
\end{split}
\end{equation}
Thus the normal form of Theorem~\ref{fin:main} gives explicit
coordinates for first-order outer automorphisms, and the parameter
bracket gives the first mixed coefficient of their commutator.
No finite-dimensional algebraic group structure on the full
automorphism group is required.

\subsection{Formal lifting over an arbitrary base ring}
\label{aut:formal-section}

For each regular vertex choose an invertible element $u_v\in vAv$.
The formulas
\begin{equation}\label{aut:corner-units}
F(v)=v,\qquad F(e)=u_{s(e)}e,\qquad
F(e^*)=e^*u_{s(e)}^{-1}
\end{equation}
define an endomorphism. Indeed, for $s(e)=s(f)=v$ one has
$F(e^*)F(f)=e^*u_v^{-1}u_vf=e^*f$, and
\[
\sum_{s(e)=v}F(e)F(e^*)
=u_v\left(\sum_{s(e)=v}ee^*\right)u_v^{-1}=v.
\]
The remaining relations follow from the corner conditions. Such
an endomorphism need not be invertible: for one loop $x$, choosing
$u=x^{-1}$ gives $x\mapsto1$.

\begin{proposition}\label{aut:formal-lift}
Every $D\in\Der_\kk(A)$ is the first coefficient of a continuous
$\kk[[t]]$-algebra automorphism of $A[[t]]$:
\[
G_t=\mathrm{id}+tD+O(t^2).
\]
This holds over every commutative base ring $\kk$.
\end{proposition}
\begin{proof}
First suppose that $D=D_p$ vanishes on the vertices. In
\eqref{aut:corner-units}, take
\[
u_v=v+tp_v,\qquad
u_v^{-1}=v-tp_v+t^2p_v^2-\cdots.
\]
The relation checks hold in $A[[t]]$. The resulting homomorphism
$A\to A[[t]]$ extends continuously to a $\kk[[t]]$-linear
endomorphism $F_t$ of $A[[t]]$, whose first coefficient is $D_p$.
The operator $B=F_t-\mathrm{id}$ increases $t$-adic order. Therefore
\[
F_t^{-1}=\sum_{j\ge0}(-B)^j
\]
is a well-defined two-sided linear inverse: each coefficient
involves only finitely many terms. The inverse of a bijective
algebra homomorphism is multiplicative.

For a general $D$, put $h_0=\sum_vD(v)v\in A$.
The vertex identities give $[h_0,v]=D(v)$, so
$D-\ad_{h_0}=D_p$ for some $p$. Then
$G_t=\operatorname{Ad}_{1+th_0}\circ F_t$ has first coefficient
$D$, by \eqref{aut:inner-first-order}.
\end{proof}

The lift is a formal path; neither the identity $G_sG_t=G_{s+t}$
nor specialization at a nonzero scalar is claimed. If
$\mathbb Q\subseteq\kk$, the formal one-parameter group
\begin{equation}\label{aut:exponential}
\exp(tD)=\sum_{n\ge0}\frac{t^nD^n}{n!}
\end{equation}
is available. The iterated Leibniz formula proves multiplicativity
and the group law. When $D$ is locally nilpotent, the series is
finite on every element. The following direct construction also
works in positive characteristic.

\subsection{Polynomial path shears}
\label{aut:shears-section}

\begin{proposition}\label{aut:path-shear}
Let $e$ be an edge and let $P$ be a path of positive length with
$s(P)=s(e)$ and $r(P)=r(e)$. Suppose that $P$ contains no occurrence
of $e$. Then the formulas
\begin{equation}\label{aut:shear-formulas}
\Phi_t(v)=v,\qquad
\Phi_t(f)=f+\delta_{f,e}\,tP,\qquad
\Phi_t(f^*)=f^*-t f^*Pe^*
\end{equation}
define a $\kk[t]$-automorphism of $A\otimes_\kk\kk[t]$.
They satisfy $\Phi_s\Phi_t=\Phi_{s+t}$ and
$\Phi_t^{-1}=\Phi_{-t}$. The first coefficient is the locally
nilpotent derivation
\[
D=\cD_{Pe^*}=\cD_{s(e),Pe^*}.
\]
\end{proposition}
\begin{proof}
The first edge of $P$ differs from $e$, so $e^*P=0$. Hence
$z=Pe^*\in s(e)As(e)$ satisfies $z^2=0$. In
\eqref{aut:corner-units} take
$u_{s(e)}=s(e)+tz$, with inverse $s(e)-tz$, and take $u_v=v$
at every other regular vertex. This gives
\eqref{aut:shear-formulas} and verifies all defining relations.

Every real edge in $P$ is fixed, so $\Phi_s(P)=P$. Also
$\Phi_s(e^*)=e^*$. The group law is immediate on real edges;
on a ghost edge the only possible additional term is
$st f^*Pe^*Pe^*=0$. Thus the group law holds and
$\Phi_{-t}$ is the inverse.

Taking the first coefficient gives $D$. Since $D(P)=D(e^*)=0$,
the operator $D^2$ vanishes on every vertex, real edge, and ghost
edge. For a product of $N$ such generators, the iterated Leibniz
formula gives $D^{N+1}=0$: every term differentiates at least one
factor twice. Thus $D$ is locally nilpotent in every characteristic.
If $\mathbb Q\subseteq\kk$, its exponential agrees with $\Phi_t$
on generators and hence on the whole algebra.
\end{proof}

The formula $\Phi_t=\mathrm{id}+tD$ here describes the generators;
higher powers of $t$ can occur on their products. The shear
preserves the standard grading when $|P|=1$. Over a field, a
nonzero scalar specialization with $|P|>1$ does not preserve it.

\subsection{Independent outer tangent directions}
\label{aut:outer-section}

\begin{proposition}\label{aut:independent-tangents}
Suppose that the finite graph has no sources. For all distinct
pairs $(P,e)$ satisfying Proposition~\ref{aut:path-shear}, the classes
\[
[\cD_{Pe^*}]\in\HH^1_\kk(A)_{|P|-1}
\]
are $\kk$-linearly independent. Any formal automorphism congruent
to the identity modulo $t$ whose first coefficient has a nonzero
finite linear combination of these classes is not an inner
automorphism of $A[[t]]$.
\end{proposition}
\begin{proof}
It suffices first to consider degree $n=|P|-1\ge0$. Suppose that
\[
p=\sum_{P,e}c_{P,e}\,Pe^*,
\]
with each summand placed in the corner $s(e)As(e)$, defines an
inner derivation. By Theorem~\ref{gen:inner}, $p=\partial h$ for
a diagonal $h$. Taking its degree-$n$ component, assume that $h$
is homogeneous of degree $n$.

If the normal expansion of $h$ contains a mixed word
$z=\alpha\beta^*$ of maximal depth
$m=\min(|\alpha|,|\beta|)\ge1$ with nonzero coefficient $c$,
choose an edge $f$ entering $u=s(\alpha)=s(\beta)$.
Such an edge exists because the graph has no sources.
The word $fzf^*$ is normal of depth $m+1$, since adjoining the
same first edge does not change the condition on the last edges.
Its coefficient in $\partial h$ is $-c$. It cannot cancel:
its common first edge determines $f$ and then $z$ uniquely;
the term $+h$ has depth at most $m$; and images of pure words
have depth at most one. On the other hand, every word $Pe^*$ in
$p$ is normal of depth one, because $P$ does not end in $e$.
This is a contradiction.

Thus $h$ consists of real closed paths if $n>0$, and of vertices
if $n=0$. The low-depth coordinates $(Q,f)$ are independent:
from $\sum_f x_ff^*=0$, right multiplication by a fixed edge
$f$ gives $x_f=0$, after which independence of real paths applies.
For a real closed path $c$, the inner relation in these coordinates is
\[
\sum_{s(f)=s(c)}(cf,f)
-\sum_{r(f)=s(c)}(fc,f).
\]
A contribution to $(P,e)$ would require $P=ce$ or $P=ec$,
both excluded by the condition on $P$. In degree zero, the
relations from vertices involve only coordinates $(e,e)$.
Hence all $c_{P,e}$ vanish. Independence across different degrees
follows from the grading.

Finally, if $G_t=\operatorname{Ad}_{U(t)}$ and $G_0=\mathrm{id}$,
then $U(0)$ is a central unit of $A$. Multiplying $U(t)$ by
$U(0)^{-1}$ does not change conjugation and allows us to write
$U(t)=1+th+O(t^2)$. Formula~\eqref{aut:inner-first-order} then
shows that the first coefficient of $G_t$ is $\ad_h$, contrary
to the assumed nonzero outer class.
\end{proof}

For linear independence alone, it is enough that $P$ neither
starts nor ends with $e$. Excluding every occurrence of $e$ is
used for the group law of the shear. The no-source hypothesis
is essential: the acyclic graph with edges $e:u\to w$,
$f:u\to v$, and $g:v\to w$ admits the shear
$e\mapsto e+tfg$, but its $\HH^1$ vanishes. For graphs with
sources, one must first apply source removal and reduce all
inner relations. A formal automorphism being outer does not
by itself imply that every nonzero scalar specialization is outer.

\subsection{Bouquets and comparison with known automorphisms}
\label{aut:bouquet}

In a bouquet of $q\ge2$ loops, distinguish a loop $a$ and let
$B_+$ be the span of nonempty words in the other $q-1$ loops.
For $F\in B_+$, the formulas
\begin{equation}\label{aut:bouquet-group}
a\longmapsto a+F,\qquad
f^*\longmapsto f^*-f^*Fa^*,
\end{equation}
with all other real edges fixed, give an embedding
$(B_+,+)\hookrightarrow\operatorname{Aut}_\kk(A)$.
Indeed, the proof of Proposition~\ref{aut:path-shear} applies
to $F$: $a^*F=0$, every real word in $F$ is fixed, and the
image of $a$ proves injectivity. Proposition~\ref{aut:independent-tangents}
also gives an injection $F\mapsto[\cD_{v,Fa^*}]$; the words
of length $m$ provide $(q-1)^m$ independent outer directions
of degree $m-1$.

Over a field, the automorphisms
$e_2\mapsto e_2+e_1p$ and $e_1^*\mapsto e_1^*-pe_2^*$
for parallel edges occur in \cite[Corollary~2.5]{KNauto},
with the bouquet case in \cite[Corollary~2.8]{KNauto}.
The associated additive-group embedding is given in
\cite[Corollary~2.9]{NSVauto}. Decomposing
$F=\sum_{f\ne a}fQ_f$ by its first edge expresses
\eqref{aut:bouquet-group} as a composition of these known
shears: each $Q_f$ uses only fixed real loops and contains no $a$.
We therefore make no separate novelty claim for the existence
of this bouquet family.

For a general path $P=fQ$, the edges $f$ and $e$ need not have
the same range, so the parallel-edge hypothesis of
\cite[Corollary~2.5]{KNauto} may fail. The direct proof above
also works over an arbitrary commutative base ring.
The cohomological conclusion used here is the independence of
the outer tangent classes, established separately in
Proposition~\ref{aut:independent-tangents}.

\section{The boundary between algebraic and \texorpdfstring{$C^*$}{C*}-derivations}
\label{analytic:boundary}

This section is over $\mathbb C$. It corrects the unrestricted
extension assertion in \cite[Theorem~4.4, p.~80]{Lopatkin}.
An algebraic derivation on a
dense subalgebra of a $C^*$-algebra need not extend to a derivation
defined on the entire $C^*$-algebra. The distinction is already
visible for one loop:
\[
\mathbb C[t,t^{-1}]\subset C(S^1).
\]
The Euler derivation $L_0=t\,d/dt$ satisfies
\[
\|L_0(t^n)\|_\infty=n,\qquad
\|t^n\|_\infty=1\qquad(n\ge1),
\]
and therefore has no bounded extension to $C(S^1)$.

In fact every complex-linear derivation
$D:C(S^1)\to C(S^1)$ is zero, with no continuity assumption.
Fix $z\in S^1$ and $f\in C(S^1)$. Put $g=f-f(z)$ and factor
$g=ab$, where
\[
a=|g|^{1/2},\qquad
b(x)=\begin{cases}
g(x)/|g(x)|^{1/2},&g(x)\ne0,\\
0,&g(x)=0.
\end{cases}
\]
Both functions are continuous: in particular, $|b|=|g|^{1/2}$
controls continuity at the zeros of $g$. Moreover $a(z)=b(z)=0$.
Since $D$ vanishes on constants, the Leibniz rule gives
\[
D(f)(z)=D(g)(z)=D(a)(z)b(z)+a(z)D(b)(z)=0.
\]
As $f$ and $z$ were arbitrary, $D=0$. Hence the nonzero Euler
derivation has no everywhere-defined extension of this kind,
even a discontinuous one.

By contrast, if $A$ is a dense subalgebra of a $C^*$-algebra $B$
and $h\in A$, then $\ad_h$ extends to $B$ and
$\|\ad_h\|\le2\|h\|$. Densely defined unbounded derivations,
their closability, and their generation of automorphism groups
require additional analytic hypotheses. The algebraic results
of this paper do not assert such properties.

\section*{Acknowledgement of AI assistance}
ChatGPT was used during the preparation of this manuscript to assist
with proof exploration, checking calculations and arguments, organizing
the exposition, and English-language editing.


\begin{thebibliography}{99}
\raggedright

\bibitem{AAJZ}
A. Alahmadi, H. Alsulami, S. K. Jain, and E. Zelmanov,
\emph{Leavitt path algebras of finite Gelfand--Kirillov dimension},
J. Algebra Appl. \textbf{11} (2012), no.~6, 1250225, 6 pp.
\href{https://doi.org/10.1142/S0219498812502258}
{doi:10.1142/S0219498812502258};
\href{https://arxiv.org/abs/1204.5258}{arXiv:1204.5258}.

\bibitem{Bez}
O. Bezushchak,
\emph{Derivations of rings of infinite matrices},
Researches in Mathematics \textbf{31} (2023), no.~2, 14--18.
\href{https://doi.org/10.15421/242310}{doi:10.15421/242310};
\href{https://arxiv.org/abs/2310.11509}{arXiv:2310.11509}.

\bibitem{CLWshort}
X.-W. Chen, H. Li, and Z. Wang,
\emph{An explicit projective bimodule resolution of a Leavitt path algebra},
preprint (2017),
\href{https://arxiv.org/abs/1711.01876}{arXiv:1711.01876}.

\bibitem{Gerritzen}
L. Gerritzen,
\emph{Modules over the algebra of the noncommutative equation \(yx=1\)},
Arch. Math. \textbf{75} (2000), 98--112.
\href{https://doi.org/10.1007/PL00000437}{doi:10.1007/PL00000437}.

\bibitem{HHToep}
Thao Nhi Nguyen Huynh and Viet Khanh Huynh,
\emph{Lie ideals and Lie derivations of the algebraic Toeplitz algebra},
preprint (2026),
\href{https://arxiv.org/abs/2609.04745}{arXiv:2609.04745}.

\bibitem{KNauto}
S. Kuroda and T. G. Nam,
\emph{Anick type automorphisms and new irreducible representations
of Leavitt path algebras},
J. Noncommut. Geom. \textbf{17} (2023), no.~3, 811--834.
\href{https://doi.org/10.4171/JNCG/489}{doi:10.4171/JNCG/489}.

\bibitem{Lopatkin}
V. Lopatkin,
\emph{Derivations of Leavitt path algebras},
J. Algebra \textbf{520} (2019), 59--89.
\href{https://doi.org/10.1016/j.jalgebra.2018.11.011}
{doi:10.1016/j.jalgebra.2018.11.011};
\href{https://arxiv.org/abs/1509.05075}{arXiv:1509.05075}.

\bibitem{LN}
V. Lopatkin and T. G. Nam,
\emph{On the homological dimensions of Leavitt path algebras with
coefficients in commutative rings},
J. Algebra \textbf{481} (2017), 273--292.
\href{https://doi.org/10.1016/j.jalgebra.2017.02.027}
{doi:10.1016/j.jalgebra.2017.02.027}.

\bibitem{NSVauto}
T. G. Nam, A. K. Srivastava, and N. T. Vien,
\emph{Automorphisms of Leavitt path algebras: Zhang twist and irreducible
representations},
J. Algebra \textbf{654} (2024), 189--234.
\href{https://doi.org/10.1016/j.jalgebra.2024.05.014}
{doi:10.1016/j.jalgebra.2024.05.014};
\href{https://arxiv.org/abs/2311.15910}{arXiv:2311.15910}.

\end{thebibliography}
\end{document}